\documentclass[a4]{article}
\usepackage{fullpage}
\usepackage{amsfonts}
\usepackage{amsmath}
\usepackage{amssymb}
\usepackage{amsthm}
\usepackage{graphicx}
\usepackage{xcolor}
\usepackage[mathscr]{eucal}

\numberwithin{equation}{section}

\title{On difference Riccati equation and continued fractions}
\author{Alexey V. Ivanov}
\date{}

\DeclareMathOperator*{\contfrac}{K}
\DeclareMathOperator*{\vprod}{\overleftarrow{\prod}}

\begin{document}
\renewcommand{\theequation}{\arabic{section}.\arabic{equation}}
\maketitle

\begin{abstract}
We study a difference Riccati equation $\Phi(x) + \rho(x)/\Phi(x-\omega) = v(x)$ with 
$1-$periodic continuos coefficients. Using continued fraction theory we investigate a problem of existence of continuos solutions for this equation. It is shown that convergence of a continued fraction representing a solution of the Riccati equation can be expressed in terms of hyperbolicity of a cocycle naturally associated to this continued fraction. We apply the critical set method to establish the uniform hyperbolicity of the cocycle and to obtain sufficient conditions for the convergence of a continued fraction giving a representation for a solution of the Riccati equation.
\end{abstract}

Keywords: difference Riccati equation, continued fraction, linear cocycle, hyperbolicity, Lyapunov exponent, critical set

MSC 2010: 37C55, 37D25, 37C40

\section{Introduction}

We study a difference Riccati equation
\begin{equation}\label{eq_1}
\Phi(x) + \rho(x)\frac{1}{\Phi(x-\omega)} = v(x),
\end{equation}
where $\rho$ and $v$ are known continuous $1$-periodic functions and $\omega\in \mathbb{T}^{1} =\mathbb{R}/\mathbb{Z}$ is irrational.

Although such equations appear in different areas of physics (e.g. control theory, theory of filters) the main interest to $(\ref{eq_1})$ comes from its relation to a linear difference equation
\begin{equation}\label{eq_2}
\vec\Psi(x+\omega) = H(x) \vec\Psi(x),\quad \vec\Psi(x)\in \mathbb{R}^{2},
\end{equation}
where $H$ is a $1$-periodic $SL(2,\mathbb{R})$-valued functions (see  e.g. \cite{BuFe94} - \cite{ABD}).

Indeed, if
\begin{equation*}
H(x) = \left(
\begin{array}{cc}
h_{11}(x) & h_{12}(x)\\
h_{21}(x) & h_{22}(x)
\end{array}
\right),
\end{equation*}
then $\Phi(x) = \frac{\Psi_{1}(x+\omega)}{\Psi_{1}(x)}$ satisfies $(\ref{eq_1})$ with
\begin{equation*}
\rho(x) = \frac{h_{12}(x)}{h_{12}(x-\omega)},\quad v(x) = h_{11}(x) + \rho(x) h_{22}(x).
\end{equation*}
Note that inverse passage from a solution of the Riccati equation to a solution of equation 
$(\ref{eq_2})$ is not simple, since in this case one needs to solve a "homological" equation
\begin{equation*}
\psi(x+\omega) - \psi(x) = \phi(x),
\end{equation*}
where $\psi(x) = \ln \Psi_{1}(x), \phi(x) = \ln \Phi(x)$. Existence of a solution for the homological equation is strongly depends on arithmetical properties of the parameter 
$\omega$ (see \cite{Russ}, \cite{Avi} for details). 

We rewrite $(\ref{eq_1})$ in a form
\begin{equation}\label{eq_3}
\Phi(x) = v(x) - \frac{\rho(x)}{\Phi(x-\omega)}.
\end{equation}
Then, representing $\Phi(x-\omega)$ by means of $(\ref{eq_3})$, one obtains
\begin{equation*} 
\Phi(x) = v(x) - \frac{\rho(x)}{v(x-\omega) - \frac{\rho(x-\omega)}{\Phi(x-2\omega)}}.
\end{equation*}
Repetition of this procedure $(n-1)$ times yields
\begin{equation*} 
\Phi(x) = b_{0}(x) + \cfrac{-1}{b_{1}(x) +  
\ldots +
\cfrac{-1}{b_{n-1}(x) + \cfrac{r_{n}(x)}{\Phi(x-n\omega)}}},
\end{equation*}
where $b_{0}(x) = v(x)$, $b_{j}(x) = \cfrac{v(x-(j+1)\omega)}{\rho(x-j\omega)}$ for $j\in \mathbb{N}$ and $r_{n}(x)=\rho(x-(n-1)\omega)/\rho(x-(n-2)\omega)$.

Thus, a formal solution of $(\ref{eq_1})$ can be represented as a continued fraction
\begin{equation*}
\Phi(x) = b_{0}(x) + \cfrac{-1}{b_{1}(x) + \cfrac{-1}{b_{2}(x) + \ldots}}  
= b_{0}(x) + \contfrac_{j=1}^{\infty}\cfrac{-1}{b_{j}(x)}.
\end{equation*}
Here we use Gauss notation, $K$, for the continued fraction.

A natural question arises: whether this continued fraction converges or not. More precisely, one may ask whether a sequence $\{\Phi_{n}\}_{n=1}^{\infty}$ with
\begin{equation*}
\Phi_{n}(x) = b_{0}(x) + \cfrac{-1}{b_{1}(x) + \cfrac{-1}{b_{2}(x) + \cdots+\cfrac{-1}{b_{n}(x)}}} = 
b_{0}(x) + \contfrac_{j=1}^{n}\cfrac{-1}{b_{j}(x)}.
\end{equation*}
converges and in what class of functions.

There is broad literature on continued fractions theory. Started from works by Euler, Lagrange, Galois this theory attracted attention of many mathematicians and found applications in various branches of mathematics such as number theory, theory of analytic functions and others. We refer the reader to classical books on this theory \cite{Per}, \cite{Khi}, \cite{JoTh} and references therein. Besides, the history of continued fractions can be found in \cite{Brez}.

Although theory of continued fractions has long history a problem of establishing the convergence of a functional continued fraction is far from its complete solution. It appears that the most general results in this direction have been obtained in a period from the middle of 19th centure to the beginning of 20th century and are usually related to such mathematicians as Worpitsky, Pringsheim, van Vleck, Stieltjes and others. In this paper we study this problem by using a well-known connection between continued fraction theory and theory of linear cocycles. In particular, we use hyperbolic property of cocycles associated to a continued fraction to establish its convergence. 

The paper is organized as follows. In section 2 we discuss general notions of continued fraction theory and its relation to one-dimensional dynamics, in particular, to theory of hyperbolic $SL(2,\mathbb{R})$-cocycles over irrational rotation. Section 3 is focused on a representation of cocycles associated to continued fractions, which simplifies further analysis of their dynamics.  In section 4 we obtain sufficient conditions for uniform hyperbolicity of cocycles.  In section 5 we apply obtained results and the critical set method to some special class of functional continued fractions and establish sufficient conditions for their convergence. 

\section{Continued fractions}

In this section we collect some known facts from the theory of continued fractions and present it's relationship with the theory of hyperbolic cocycles. Our main references on the theory of continued fractions are \cite{Khi}, \cite{JoTh}. 

\newtheorem{defs}{Definition}
\begin{defs}
 A sequence $\{(a_{j}, b_{j}), f_{j}\}$, where $b_{0}, a_{j}, b_{j}, f_{j}\in \mathbb{R}\,\, \forall j\in \mathbb{N}$ is called a continued fraction if 
\begin{equation*}
f_{n} = b_{0} +  \cfrac{a_{1}}{b_{1} + \cfrac{a_{2}}{b_{2} + \ldots+\cfrac{a_{n}}{b_{n}}}} = b_{0} + \contfrac_{j=1}^{n}\cfrac{a_{j}}{b_{j}},\quad
\forall\, n\in\mathbb{N}.
\end{equation*}
The number $f_{n}$ is called the $n-$th convergent of the continued fraction.
\newline
One says that a continued fraction $\{(a_{j}, b_{j}), f_{j}\}$ converges if there exists finite or infinite limit
\begin{equation*}
\lim\limits_{n\to \infty}f_{n}=f_{*}.
\end{equation*}
\end{defs}
It may occur that for a two different continued fractions $\{(a_{j}, b_{j}), f_{j}\}$, $\{(a'_{j}, b'_{j}), f'_{j}\}$ the corresponding convergents coincide, i.e.
\begin{equation*}
f_{j} = f'_{j},\,\, \forall j\in \mathbb{N}.
\end{equation*}
In this case one says that the continued fractions $\{(a_{j}, b_{j}), f_{j}\}$, $\{(a'_{j}, b'_{j}), f'_{j}\}$ are equivalent. In particular, if for a given continued fraction $\{(a_{j}, b_{j}), f_{j}\}$ such that $a_{j}\neq 0,\,\,\forall j\in\mathbb{N}$ one sets
\begin{equation*}
a'_{j} = -1,\quad b'_{j}=r_{j}b_{j},
\end{equation*}
where $r_{0}=1$ and
\begin{equation*}
r_{2n}=\frac{\prod\limits_{s=1}^{n}a_{2s-1}}{\prod\limits_{s=1}^{n}a_{2s}},\quad
r_{2n+1}=-\frac{\prod\limits_{s=1}^{n}a_{2s}}{\prod\limits_{s=0}^{n}a_{2s+1}},
\end{equation*}
then we obtain an equivalent continued fraction \cite{JoTh}.

Another remarkable fact is the existence of inductive procedure, which enables constructing the $n-$th convergents \cite{Khi}. Namely, for a continued fraction 
$\{(a_{j}, b_{j}), f_{j}\}$ define 
\begin{align}\label{f_AB}
\nonumber
&p_{n} = b_{n}p_{n-1} + a_{n}p_{n-2},\\
&q_{n} = b_{n}q_{n-1} + a_{n}q_{n-2},\\
\nonumber
&p_{0}=b_{0},\,\, p_{-1} = 1,\,\, q_{0}=1,\,\, q_{-1}=0.
\end{align}
Then 
\begin{equation}\label{n_conv}
f_{n} = \frac{p_{n}}{q_{n}}.
\end{equation}

Although the development of continued fraction theory started in 16th centure, there is a lack of sufficient conditions guaranteeing convergence of a continued fraction of
general type. The most powerfull tools for establishing the convergence of continued fractions with real entries are due to the following two theorems \cite{JoTh}:
\newtheorem*{thmSS}{Theorem (L. Seidel, M. A. Stern, 1846-1848)}
\begin{thmSS}
A continued fraction $\{(a_{j}, b_{j}), f_{j}\}$ converges to a finite value if
\begin{equation*}
a_{j}=1, \quad b_{j}>0,\quad \forall j\in \mathbb{N}
\end{equation*}
and
\begin{equation*}
\sum\limits_{j=1}^{\infty}b_{j} = \infty.
\end{equation*}
\end{thmSS}

\newtheorem*{thmP}{Theorem (A. Pringsheim, 1899)}
\begin{thmP}
A continued fraction $\{(a_{j}, b_{j}), f_{j}\}$ converges to a finite value if
\begin{equation*}
\vert b_{j}\vert \ge \vert a_{j}\vert + 1,\quad \forall j\in \mathbb{N}.
\end{equation*}
\end{thmP}
As a consequence to Pringsheim's theorems one obtains Worpitsky's sufficient condition, which establishes convergence of a continued fraction under conditions
\begin{equation}\label{Worp_cond}
\vert a_{j}\vert = 1,\quad \vert b_{j}\vert \ge 2,\quad \forall j\in \mathbb{N}.
\end{equation}

One may note that formula (\ref{n_conv}) can be rewritten in terms of a cocycle. Indeed, for a given continued fraction $\{(a_{j}, b_{j}), f_{j}\}$ we consider a sequence of matrices
\begin{equation}\label{A_trans}
A_{n} = \left(
\begin{array}{cc}
b_{n} & a_{n}\\
1 & 0
\end{array}
\right).
\end{equation}
This sequence generates a cocycle 
\begin{equation*}
\quad M_{n} = A_{n}A_{n-1}\cdots A_{1}=
\vprod\limits_{k=1}^{n}A_{k},\quad n\in \mathbb{N},
\end{equation*}
where the left arrow denotes the order of multiplication.

Then (\ref{f_AB}) can be rewritten as
\begin{align*}
&\vec p_{n} = A_{n}\vec p_{n-1} = M_{n}\vec p_{0},\quad
\vec q_{n} = A_{n}\vec q_{n-1} = M_{n}\vec q_{0},\\
&\vec p_{n} = \left(
\begin{array}{c}
p_{n}\\
p_{n-1}
\end{array}
\right),\quad
\vec q_{n} = \left(
\begin{array}{c}
q_{n}\\
q_{n-1}
\end{array}
\right).
\end{align*}
Hence, the $n$-th convergent reads
$$
f_{n} = \frac{(M_{n}\vec e_{2}, \vec e_{1})}{(M_{n}\vec e_{1}, \vec e_{1})},
$$
where $\{\vec e_{1}, \vec e_{2}\}$ is the standard basis in $\mathbb{R}^{2}$.

It is important for further analysis that transformation $A_{n}\in SL(2,\mathbb{R})$ provided $a_{n} = -1$.
\newline
Now we consider functional continued fractions $\{(a_{j}, b_{j}), f_{j}\}$, whose elements $a_{j}, b_{j}$ are $1$-periodic continuous functions of a real variable $x$ or equivalently we may suppose that $a_{j}, b_{j}$ are continuous functions on a circle $\mathbb{T}^{1} = \mathbb{R}/\mathbb{Z}$.
\begin{defs}
We say that a functional continued fraction $\{(a_{j}, b_{j}), f_{j}\}$ is generated by a pair $(b, \omega)$, where $b$ is a $1$-periodic continuous function and $\omega\in (0,1)$ is an irrational number if 
\begin{equation*}
a_{j}(x)\equiv -1,\quad b_{j}(x) = b\bigl(x-(j-1)\omega\bigr),\quad \forall j\in \mathbb{N},\quad 
\forall x\in \mathbb{R}.
\end{equation*}
\end{defs}
Thus, in that case the associated cocycle $M_{n}$
\begin{equation*}
M_{n}(x) = \vprod\limits_{k=1}^{n}A_{k}(x),\quad 
A_{k}(x) = \left(
\begin{array}{cc}
b\bigl(x-(k-1)\omega\bigr) & -1\\
1 & 0
\end{array}
\right)
\end{equation*}
becomes an $SL(2,\mathbb{R})$-cocycle. 
The following theorem holds true.
\newtheorem{thms}{Theorem}
\begin{thms}
Assume that a given functional continued fraction $\{(a_{j}, b_{j}), f_{j}\}$ is generated by a pair $(b, \omega)$ with function $b\in C(\mathbb{T}^{1},\mathbb{R})$ and $\omega\in (0,1)\cap \mathbb{R}\setminus\mathbb{Q}$.
If there exist positive constants $C_{0}$ and $\lambda_{0}>1$ such that the cocycle associated with this continued fraction admits an estimate
\begin{equation}\label{coc_est}
\Vert M_{n}(x)\Vert \ge C_{0} \lambda_{0}^{n},
\end{equation}
then the continued fraction converges and
$f_{*} = \lim\limits_{n\to \infty}f_{n}\in C(\mathbb{T}^{1},\mathbb{\hat R})$, where
$\mathbb{\hat R} = \mathbb{R}\cup\{\infty\}$.
\end{thms}
{\bf PROOF:} It follows from results of J.-C. Yoccos \cite{Yoc} (see also \cite{AvBo}) that in the case of $SL(2,\mathbb{R})$-cocycles estimate (\ref{coc_est}) is equivalent to the uniform hyperbolicity of the cocycle. Then by Oseledets theorem \cite{Oseled} there exist continuous maps
\begin{equation*}
E^{s,u}:\mathbb{T}^{1}\to Gr(1,2),
\end{equation*}
where $Gr(1,2)$ stands for the set of one-dimensional subspaces of 
$\mathbb{R}^{2}$, such that 
\begin{equation*}
M_{0}(x)E^{s,u}(x)=E^{s,u}(x-\omega),\quad 
E^{s}(x) \oplus E^{u}(x) = \mathbb{R}^{2}
\end{equation*} 
and 
\begin{align}\label{Th_O}
\left\Vert M_{n}(x)\vert_{E^{s}(x)}\right\Vert \le C{\rm e}^{-\Lambda n}, \quad n\ge 0,\\
\nonumber
\left\Vert M_{n}(x)\vert_{E^{u}(x)}\right\Vert \le C{\rm e}^{-\Lambda n},\quad n\le 0
\end{align}
Let $\vec e^{s,u}(x)\in E^{s,u}(x)$ be unit vectors. Due to continuity and transversality of $E^{s,u}$  we may assume that $\vec e^{s,u}$ depend continuously on $x\in \mathbb{T}^{1}$. Then
\begin{equation*}
\vec e_{k} = c_{k}(x)\vec e^{u}(x) + d_{k}(x)\vec e^{s}(x),\quad k=1,2,
\end{equation*}
where $c_{k}, d_{k}$ are continuous functions. Then, using estimates (\ref{Th_O}), one obtains
\begin{equation*}
\lim\limits_{n\to \infty}
\frac{(M_{n}(x)\vec e_{2}, \vec e_{1})}{(M_{n}(x)\vec e_{1}, \vec e_{1})}=
\frac{c_{2}(x)}{c_{1}(x)}
\end{equation*}
Applying transversality of subspaces $E^{s,u}$, we conclude that $c_{1}, c_{2}$ cannot vanish simultaneously. This finishes the proof.$\qed$

{\bf Remark:} One may note that if $E^{s}(x)$ is transversal to $E_{1}=\{t \vec e_{1}, t\in \mathbb{R}\},\,\, \forall x\in \mathbb{T}^{1}$ then $c_{1}$ does not vanish and $f_{*}\in C(\mathbb{T}^{1},\mathbb{R})$.

Thus, Theorem 1 says that uniform hyperbolicity of the cocycle $M_{n}$ implies convergence of the functional continued fraction.

\section{Decomposition of cocycles}

In this section we describe a procedure, which allows passage from a cocycle associated to a continued fraction to a cocycle which possesses the same hyperbolic properties, but has a special representation convenient for analysis.  We will suppose in what follows that a continued fraction $\{(a_{j}, b_{j}), f_{j}\}$ satisfies 
\begin{equation*}
a_{j}\neq 0\quad \forall\, j\in \mathbb{N}.
\end{equation*}
Hence, passing to an equivalent continued fraction (if necessary) a transformation (\ref{A_trans}) can be written in the form
\begin{equation*}
A_{n} = \left(
\begin{array}{cc}
b_{n} & -1\\
1 & 0
\end{array}
\right)\in SL(2,\mathbb{R}).
\end{equation*}

To analyze behavior of $M_{n}$ for large $n$ we represent the transformation $A_{n}$ by using of singular value decomposition (SVD). It can be done in different ways, however, we choose the following SVD for $A_{n}$:
\begin{equation*}
A_{n} = \textrm{sign}(b_{n}) R(\theta_{n}) Z(\mu_{n}) R(\theta_{n}), 
\end{equation*}
where 
\begin{equation*}
R(\varphi) = \left(
\begin{array}{cc}
\cos \varphi & \sin \varphi\\
-\sin \varphi & \cos \varphi
\end{array}
\right), 
\,\,
Z(\lambda) = \left(
\begin{array}{cc}
\lambda & 0\\
0 & \lambda^{-1}
\end{array}
\right)
\end{equation*}
and
\begin{align*}
&\mu_{n} = \frac{1}{2}\left(\sqrt{b_{n}^{2}+4}+\vert b_{n}\vert\right),\\ 
\nonumber
&\tan \theta_{n} = \frac{\textrm{sign}(b_{n})}{2}\left(\sqrt{b_{n}^{2}+4}-\vert b_{n}\vert\right).
\end{align*}

The reason why we choose this representation is that 
\begin{equation*}
\mu_{n} \ge 1,\quad \forall n\in \mathbb{N}.
\end{equation*}
This fact simplifies analysis of cocycle dynamics. Besides, one may note that
\begin{equation*}
\vert \tan\theta_{n}\vert \le 1,\quad \forall n\in \mathbb{N}
\end{equation*}
and $\mu_{n}$ is a continuous function of $b_{n}$. However, $\theta_{n}$ is discontinuous with respect to $b_{n}$ at the point $b_{n}=0$. The latter property makes further analysis more difficult.

To overcome this obstacle we consider a sequence of matrices $A^{(2)}_{n} = A_{2n}A_{2n-1}$. Note that a cocycle, $M_{n}^{(2)}$, generated by $A^{(2)}_{n}$ is connected with the cocycle $M_{n}$ in the following way
\begin{equation*}
M_{n}^{(2)} = M_{2n}.
\end{equation*}
It is known that hyperbolic properties of both cocycles are the same. A direct applcation of singular value decomposition to $A^{(2)}_{n}$ leads to
\begin{equation*}
A^{(2)}_{n} = R(\varphi_{n}) Z(\lambda_{n}) R(\chi_{n}), 
\end{equation*}
with 
\begin{align}\label{SVD_2}
\nonumber
&\lambda_{n} = \frac{1}{2}\left(\sqrt{4+h_{n}^{2}}+\sqrt{h_{n}^{2}}\right),\quad
%\nonumber 
h_{n}^{2} = (b_{2n-1}-b_{2n})^{2}+b_{2n-1}^{2}b_{2n}^{2},\\
&\cot \varphi_{n} = \epsilon_{2n-1}\sqrt{\frac{\lambda_{n}^{2}-b_{2n-1}^{2}-1}{b_{2n-1}^{2}+1-\lambda_{n}^{-2}}},\quad
%\nonumber
\cot \chi_{n} = \epsilon_{2n}\sqrt{\frac{\lambda_{n}^{2}-b_{2n}^{2}-1}{b_{2n}^{2}+1-\lambda_{n}^{-2}}},\\
\nonumber
&\epsilon_{2n-1} = \textrm{sign}\bigl(b_{2n}(1+b_{2n-1}^{2})-b_{2n-1}\bigr),\quad
%\nonumber.
\epsilon_{2n} = \textrm{sign}\bigl(b_{2n-1}(1+b_{2n}^{2})-b_{2n}\bigr).
\end{align}

Properties of decomposition (\ref{SVD_2}) can be described by the following lemma.
\newtheorem{lems}{Lemma}
\begin{lems}
Function $\lambda_{n}$ is continuous with respect to $b_{2n-1}, b_{2n}$ and satisfies
\begin{equation*}
\lambda_{n}\ge 1,\quad \forall n\in \mathbb{N}.
\end{equation*}
If additionally the following condition
\begin{equation*}
(H_{1})\quad\quad
\inf\limits_{n}\biggl( (b_{2n-1}-b_{2n})^{2}+b_{2n-1}^{2}b_{2n}^{2}\biggr) = h_{*}^{2} > 0
\end{equation*}
is fullfilled, then functions $\varphi_{n}$, $\chi_{n}$ are continuous with respect to $b_{2n-1}, b_{2n}$ and 
there exists a positive constant $\Lambda_{0}>1$ such that
\begin{equation*}
\lambda_{n}\ge \Lambda_{0},\quad \forall n\in \mathbb{N}.
\end{equation*}
\end{lems}
{\bf PROOF:} The continuity of $\lambda_{n}$ immediately follows from formulae (\ref{SVD_2}). Also one may note that $\lambda_{n}$ increases with respect to $h_{n}$, while condition 
$(H_{1})$ implies $h_{n}\ge h_{*}$. Hence, one may set
\begin{equation*}
\Lambda_{0} = \frac{1}{2}\left(\sqrt{4+h_{*}^{2}}+\sqrt{h_{*}^{2}}\right).
\end{equation*}
To prove continuity of $\varphi_{n}$, $\chi_{n}$ we apply formula for $\lambda_{n}$  and obtain that 
\begin{equation*}
\lambda_{n}^{2} - b_{2n-1}^{2} - 1 = 
2\frac{\left(b_{2n-1}^{2} + 1\right)^{2}\left(b_{2n} - b_{2n}^{(*)}\right)^{2}}
{\sqrt{h_{n}^{2}}\sqrt{h_{n}^{2} + 4} - h_{n}^{2} + 2b_{2n-1}^{2}}, 
\end{equation*}
where we introduced 
\begin{equation*}\label{eq_b_2k}
b_{2n}^{(*)} = \frac{b_{2n-1}}{1+b_{2n-1}^{2}}.
\end{equation*}
The same substitution also yields
\begin{equation*}
b_{2n-1}^{2} + 1 - \lambda_{n}^{-2} = \frac{1}{2}
\left(\sqrt{h_{n}^{2}}\sqrt{h_{n}^{2} + 4} - h_{n}^{2} + 2b_{2n-1}^{2}\right).
\end{equation*} 
Since $b_{2n-1}, b_{2n}$ enter into $\lambda_{n}$ symmetrically, one obtains similar formulae by interchanging $b_{2n-1}$ and $b_{2n}$. Namely,
\begin{equation*}
\lambda_{n}^{2} - b_{2n}^{2} - 1 = 
2\frac{\left(b_{2n}^{2} + 1\right)^{2}\left(b_{2n-1} - b_{2n-1}^{(*)}\right)^{2}}
{\sqrt{h_{n}^{2}}\sqrt{h_{n}^{2} + 4} - h_{n}^{2} + 2b_{2n}^{2}}, 
\end{equation*}
where  
\begin{equation*}\label{eq_b_2k-1}
b_{2n-1}^{(*)} = \frac{b_{2n}}{1+b_{2n}^{2}}
\end{equation*}
and
\begin{equation*}
b_{2n}^{2} + 1 - \lambda_{n}^{-2} = \frac{1}{2}
\left(\sqrt{h_{n}^{2}}\sqrt{h_{n}^{2} + 4} - h_{n}^{2} + 2b_{2n}^{2}\right).
\end{equation*}
Notice also that 
\begin{equation*}
\epsilon_{2n-1} = {\rm sign}\left(b_{2n} - b_{2n}^{(*)}\right),\quad
\epsilon_{2n} = {\rm sign}\left(b_{2n-1} - b_{2n-1}^{(*)}\right),\quad \forall n\in \mathbb{N}
\end{equation*}
and
\begin{equation*}
\lambda_{n}^{2}-b_{2n-1}^{2}-1\ge 0,\quad
\lambda_{n}^{2}-b_{2n}^{2}-1\ge 0,\quad \forall n\in \mathbb{N}.
\end{equation*}
Taking this into account we rewrite expressions for $\cot\varphi_{n}$ and $\cot\psi_{n}$ as
\begin{equation}\label{f_varphi_psi_new}
\cot\varphi_{n} = \frac{b_{2n-1}^{2}+1}{b_{2n-1}^{2}+1-\lambda_{n}^{-2}}
\left(b_{2n} - b_{2n}^{(*)}\right),\quad
\cot\psi_{n} = \frac{b_{2n}^{2}+1}{b_{2n}^{2}+1-\lambda_{n}^{-2}}
\left(b_{2n-1} - b_{2n-1}^{(*)}\right).
\end{equation}
We conclude that continuity of $\varphi_{n}$, $\chi_{n}$ follows from 
$(\ref{f_varphi_psi_new})$. This finishes the proof. $\qed$

In the rest of the paper we will assume that condition $(H_{1})$ from Lemma 1
takes place. For example, it is satisfied if the sequence $\{b_{j}\}$ is such that 
the Worpitsky's condition (\ref{Worp_cond}) holds or if for all 
$j\in \mathbb{N}$ the following implication is valid
\begin{equation*}
\vert b_{j}\vert < 2 \Rightarrow \vert b_{j+1}\vert \ge 2.
\end{equation*}
Finally, we introduce a sequence of matrices
\begin{align}\label{coc_B}
&B_{n} = R(\Phi_{n}) Z(\lambda_{n}),\\
\nonumber
&\Phi_{n}=\psi_{n+1}+\varphi_{n},\quad n\in \mathbb{N}.
\end{align}
It is to be noted that hyperbolic properties of cocycles generated by $\{A_{n}\}$ and 
$\{B_{n}\}$ are equivalent. On the other hand, cocycles generated by matrices (\ref{coc_B}) have simpler form and were studied in many papers (e.g. \cite{Her}, \cite{LSY}). We use results of \cite{Iva21}, \cite{Iva23} to establish hyperbolicity of the cocycle generated by matrices (\ref{coc_B}).

Due to represenation (\ref{coc_B}) a cocycle generated by the sequence $\{B_{n}\}$ is a product of interchanged matrices $R$ and $Z$, i.e.
\begin{equation}\label{coc_B_2}
M_{n} = \vprod\limits_{k=1}^{n}R(\Phi_{k})Z(\lambda_{k}).
\end{equation}
Matrices $R$ and $Z$ possesses an obvious property:
\begin{equation*}
R(x_{1}+x_{2}) = R(x_{1})R(x_{2}),\quad
Z(x_{1}+x_{2}) = Z(x_{1})Z(x_{2}),\quad
\forall\, x_{1}, x_{2}\in \mathbb{R}
\end{equation*}
and, hence, matrices $R$ (respectively, $Z$) commute. 
On the other hand, $R(x_{1})Z(x_{2})\neq Z(x_{2})R(x_{1})$ (in general). However, we have the following lemma
\begin{lems}
A product $P_{1} = Z(\lambda_{2}) R(\varphi) Z(\lambda_{1})$ can be represented as
\begin{equation*}
P_{1} = R(\psi) Z(\mu) R(\chi),
\end{equation*}
where
\begin{align}\label{f_P1}
\nonumber
&\mu = S(\lambda_{1}, \lambda_{2}, \varphi) = \frac{1}{2}\left(m(\lambda_{1}, \lambda_{2}, \varphi) + \sqrt{m^{2}(\lambda_{1}, \lambda_{2}, \varphi)-4}\right),\\
&m(\lambda_{1}, \lambda_{2}, \varphi) = (\lambda_{1}\lambda_{2} + \lambda_{1}^{-1}\lambda_{2}^{-1})
\sqrt{\cos^{2}\varphi + \beta^{2}\sin^{2}\varphi},\quad
\beta = 
\frac{\lambda_{1}^{-2}+\lambda_{2}^{-2}}{1+\lambda_{1}^{-2}\lambda_{2}^{-2}},\\
\nonumber
&\tan \chi = T(\lambda_{1}, \lambda_{2}, \varphi),\quad
\tan \psi = T(\lambda_{2}, \lambda_{1}, \varphi),\quad
T(\lambda_{1}, \lambda_{2}, \varphi) = 
\epsilon(\varphi)\sqrt{1+z^{2}(\lambda_{1}, \lambda_{2}, \varphi)} -
z(\lambda_{1}, \lambda_{2}, \varphi),\\
\nonumber
&z(\lambda_{1}, \lambda_{2}, \varphi) = \frac{\lambda_{1}^{2}}{2(1-\lambda_{2}^{-4})}
\left((1-\lambda_{1}^{-4}\lambda_{2}^{-4})\cot\varphi + 
(\lambda_{2}^{-4}-\lambda_{1}^{-4})\tan\varphi\right),\quad
\epsilon(\varphi) = {\rm sign}(2\varphi).
\end{align}
\end{lems}
{\bf PROOF:} Formulae (\ref{f_P1}) can be verified by their direct substitution into $P_{1}$.  $\qed$
\newline
We remark here that Lemma 2 is analogous to Lemma 1 from \cite{Iva23}. However, formulae (\ref{f_P1}) seem to be more convenient for analysis. In particular, as an immediate consequence, one has:
\begin{align*}
&S(\lambda_{1}, \lambda_{2}, \varphi) = 
S(\lambda_{2}, \lambda_{1}, \varphi) = 
S(\lambda_{1}^{-1}, \lambda_{2}^{-1}, \varphi) = 
S(\lambda_{1}, \lambda_{2}, -\varphi),\\
&T(\lambda_{1}, \lambda_{2}, \varphi) = - T(\lambda_{1}, \lambda_{2}, -\varphi).
\end{align*}

\section{Hyperbolicity and convergence}

In this section we obtain sufficient conditions for the convergence of continued fractions by establishing hyperbolicity of the correponding cocycles.

First, we apply Lemma 2 to a cocycle (\ref{coc_B_2}) to get the following 
\begin{lems}
Assume $\Lambda_{0}, C_{\lambda}, \delta$ be positive constants such that
\begin{equation}\label{cond_H2_1}
\Lambda_{0} > 1,\quad C_{\lambda} > \frac{\Lambda_{0}}{\sqrt{\Lambda_{0}^{2}-1}},\quad
\delta > \frac{\Lambda_{0}^{2}}{\Lambda_{0}^{2}-1}
\end{equation}
and 
\begin{equation}\label{cond_H2_2}
C_{\lambda} \ge \frac{\delta - 1 + \Lambda_{0}^{2}}{(\delta-1)\Lambda_{0}^{2} - \delta}.
\end{equation}
If parameters of a cocycle (\ref{coc_B_2}) satisfy
\begin{equation}\label{cond_H2_3}
\lambda_{n}\ge \Lambda_{0},\quad 
\Lambda_{0}\bigl\vert\cot(\Phi_{n})\bigr\vert \ge \delta C_{\lambda},\quad
\forall\, n\in \mathbb{N},
\end{equation}
then the cocycle admits an estimate
\begin{equation*}
\Vert M_{n}\Vert \ge \Lambda_{0} \hat C_{\lambda}^{n-1},\quad \hat C_{\lambda} = \frac{\Lambda_{0} C_{\lambda}}{\sqrt{\Lambda_{0}^{2} + C_{\lambda}^{2}}},\quad
\forall\, n\in \mathbb{N}
\end{equation*}
\end{lems}
{\bf PROOF:} First we apply Lemma 2 to the product  (\ref{coc_B_2}) and represent it as
\begin{equation*}
M_{n} = R(\Phi_{n}+\psi_{n})Z(\mu_{n})R(\chi_{n}), 
\end{equation*}
where 
\begin{equation}\label{step_1}
\mu_{1} = \lambda_{1},\quad \psi_{1} = 0,\quad \chi_{1} = 0
\end{equation} 
and 
\begin{equation}\label{step_n}
\mu_{n} = S(\mu_{n-1}, \lambda_{n}, \Phi_{n-1}+\psi_{n-1}),\quad
\tan\psi_{n} = T(\lambda_{n}, \mu_{n-1}, \Phi_{n-1}+\psi_{n-1}).
\end{equation}
Note that to prove Lemma 3 it suffices to prove the following estimates:
\begin{equation}\label{f_ind}
\mu_{n}\ge \Lambda_{0}\hat C_{\lambda}^{n-1},\quad \Lambda_{0}\vert\cot(\Phi_{n}+\psi_{n})\vert \ge C_{\lambda},\quad \forall n\in\mathbb{N}.
\end{equation}
We argue by induction. For $n=1$ estimates (\ref{f_ind}) follow from (\ref{step_1}), (\ref{cond_H2_1}) and (\ref{cond_H2_3}). Assume that (\ref{f_ind}) are valid for $n = k$.
Note that $1 < \hat C_{\lambda} < C_{\lambda}$ due to (\ref{cond_H2_1}). Besides, taking into account Lemma 2, one has
\begin{equation*}
\mu_{k+1}\ge \lambda_{k+1}\mu_{k}\bigl\vert \cos(\Phi_{k}+\psi_{k})\bigr\vert.
\end{equation*}
Hence, the first estimate of (\ref{cond_H2_3}) together with estimates (\ref{f_ind}) for $n = k$ imply
\begin{equation*}
\mu_{k+1}\ge \Lambda_{0}\hat C_{\lambda}^{k}.
\end{equation*}
To prove the second estimate (\ref{f_ind}) we introduce
\begin{equation*}
z_{k+1} = z(\lambda_{k+1}, \mu_{k}, \Phi_{k}+\psi_{k}).
\end{equation*}
Then one gets
\begin{equation*}
z_{k+1} = \frac{1}{2(1-\mu_{k}^{-4})}\left(\left(1-\lambda_{k+1}^{-4}\mu_{k}^{-4}\right) u_{k} - \left(1-\frac{\lambda_{k+1}^{4}}{\mu_{k}^{4}}\right)u_{k}^{-1}\right),\quad
u_{k} = \lambda_{k+1}^{2}\cot \Phi_{k}.
\end{equation*}
Using formulae (\ref{f_P1}) and (\ref{f_ind}) for $n = k$, we obtain
\begin{equation}\label{f_psi_k+1}
\bigl\vert\tan \psi_{k+1}\bigr\vert < \frac{1}{2\vert z_{k+1}\vert} < \frac{1-\mu_{k}^{-4}}{1-\lambda_{k+1}^{-4}\mu_{k}^{-4}}\cdot \frac{1}{u_{k}-u_{k}^{-1}} < 
\frac{\Lambda_{0} C_{\lambda}}{\Lambda_{0}^{2} C_{\lambda}^{2} - 1}.
\end{equation}
Note that by assumption (\ref{cond_H2_3}) $\vert\tan\Phi_{k+1}\vert \le \frac{\Lambda_{0}}{\delta C_{\lambda}}$. This estimate together with (\ref{f_psi_k+1}) and (\ref{cond_H2_3}) yield
\begin{align*}
&\bigl\vert \tan\Phi_{k+1}\cdot \tan\psi_{k+1}\bigr\vert \le \frac{1}{\delta}\cdot \frac{\Lambda_{0}^{2}}{\Lambda_{0}^{2}C_{\lambda}^{2}-1} < 1,\\
&\bigl\vert \tan\Phi_{k+1}\bigr\vert + \bigl\vert\tan\psi_{k+1}\bigr\vert \le \frac{\Lambda_{0}}{\delta C_{\lambda}} + \frac{\Lambda_{0}C_{\lambda}}{\Lambda_{0}^{2}C_{\lambda}^{2}-1}
\end{align*}
and, consequently,
\begin{equation}\label{f_Phi+psi_k+1}
\left\vert\tan\bigl(\Phi_{k+1}+\psi_{k+1}\bigr)\right\vert \le \frac{\frac{\Lambda_{0}}{\delta C_{\lambda}} + \frac{\Lambda_{0}C_{\lambda}}{\Lambda_{0}^{2}C_{\lambda}^{2}-1}}
{1-\frac{1}{\delta}\cdot \frac{\Lambda_{0}^{2}}{\Lambda_{0}^{2}C_{\lambda}^{2}-1}}.
\end{equation}
We apply inequality (\ref{cond_H2_2}) to estimate the right hand side of (\ref{f_Phi+psi_k+1}) and obtain
\begin{equation*}
\left\vert\tan\bigl(\Phi_{k+1}+\psi_{k+1}\bigr)\right\vert \le \frac{\frac{\Lambda_{0}}{\delta C_{\lambda}} + \frac{\Lambda_{0}C_{\lambda}}{\Lambda_{0}^{2}C_{\lambda}^{2}-1}}
{1-\frac{1}{\delta}\cdot \frac{\Lambda_{0}^{2}}{\Lambda_{0}^{2}C_{\lambda}^{2}-1}}\le \frac{\Lambda_{0}}{C_{\lambda}}.
\end{equation*}
This finishes the proof.$\qed$
\newline
We remark that Lemma 3 was proved in \cite{Iva23} under assumption $\Lambda_{0}\gg 1$. In this case conditions (\ref{cond_H2_1}) - (\ref{cond_H2_3}) are fullfilled if
\begin{equation*}
\lambda_{n}\ge \Lambda_{0},\quad
\bigl\vert \cot(\Phi_{n})\bigr\vert \ge \frac{C_{0}}{\Lambda_{0}},\quad
\forall\, n\in\mathbb{N}
\end{equation*}
for some positive constant $C_{0}$.

As a direct consequence of Lemma 3 one obtains
\begin{thms}
Assume that for a given continued fraction $\{(-1, b_{j}), f_{j}\}$ the condition $(H_{1})$ holds true. Let $\Lambda_{0}$ be a constant guaranteed by Lemma 1 and $\{B_{k}\}$ be the sequence of matrices (\ref{coc_B}) associated to the continued fraction.
If there exist positive constants $C_{\lambda}$ and $\delta$ satisfying (\ref{cond_H2_1}), (\ref{cond_H2_2}) such that 
\begin{equation*}
\Lambda_{0}\bigl\vert \cot(\Phi_{n})\bigr\vert \ge \delta C_{\lambda},\quad \forall\, n\in\mathbb{N},\quad\quad (H_{2})
\end{equation*}
then the cocycle generated by $\{B_{k}\}$ is hyperbolic and the continued fraction converges.
\end{thms}
It has to be noted that conditions of Pringsheim's theorem imply conditions 
$(H_{1,2})$, so Theorem 2 provides a more general sufficient condition for the convergence of a continued fraction. However, if we consider the case $(-1)^{j}b_{j}>0$ what corresponds to the Seidel-Stern theorem, it is not clear which of the conditions is stronger.

We apply this theorem to the case of functional continued fractions. If a functional continued fraction is generated by a pair $(b, \omega)$ with $b\in C(\mathbb{T}^{1}, \mathbb{R})$ and $\omega\in (0,1)\setminus \mathbb{Q}$, then trajectory of any point $x\in \mathbb{T}^{1}$ under rotation $x\mapsto x-\omega$  is dense and condition $(H_{1})$ can be replaced by
\begin{equation*}
(H'_{1})\quad\quad
\min\limits_{x\in \mathbb{T}^{1}} 
\biggl( \bigl(b(x)-b(x-\omega)\bigr)^{2}+b^{2}(x) b^{2}(x-\omega)\biggr) = 
h_{*}^{2} > 0.
\end{equation*}
In that case (similarly to Lemma 1) there exists a constant $\Lambda_{0} > 1$ such that
\begin{equation*}
\lambda_{n}(x)\ge \Lambda_{0},\quad \forall\, n\in\mathbb{N},\,\, 
\forall\, x\in \mathbb{T}^{1}.
\end{equation*}
Thus, we arrive at the following
\newtheorem{cors}{Corollary}
\begin{cors}
Let a functional continued fraction $\{(a_{j}, b_{j}), f_{j}\}$ be generated by a pair $(b, \omega)$ with $b\in C(\mathbb{T}^{1}, \mathbb{R})$ and $\omega\in (0,1)\setminus \mathbb{Q}$.
If, additionally, condition $(H'_{1})$ is fullfilled and there exist positive constants $C_{\lambda}$ and $\delta$ satisfying (\ref{cond_H2_1}), (\ref{cond_H2_2}) such that 
\begin{equation*}
\Lambda_{0} \min\limits_{x\in \mathbb{T}^{1}}\bigl\vert \cot(\Phi_{1}(x))\bigr\vert \ge \delta C_{\lambda},\quad\quad (H'_{2})
\end{equation*}
then the functional continued fraction $\{(-1, b_{j}), f_{j}\}$ converges and $f_{*} = \lim\limits_{n\to \infty}f_{n}\in C(\mathbb{T}^{1},\mathbb{\hat R})$.
\end{cors}
{\bf PROOF:} One may note that under condition of the Corollary hypotheseses 
$(H_{1,2})$ are satisfied uniformly with respect to $x$. Then application of Theorem 1 finishes the proof. $\qed$

It has to be emphasized that conditions $(H_{2})$ and $(H'_{2})$ are rather restrictive. In particular, $(H'_{2})$ is violated if a continued fraction is generated by a pair $(b, \omega)$ and 
the function $b\in C(\mathbb{T}^{1}, \mathbb{R})$ has zeroes. 

In the rest of the paper we will assume that hypothesis $(H_{2})$ is not satisfied.  To establish convergence of a continued fraction in this case we formulate two lemmae, which illustrate the main idea.

For a sequence of $SL(2,\mathbb{R})$-matrices $\{B_{n}\}_{n=1}^{\infty}$ of the form 
$(\ref{coc_B})$ and an increasing sequence of indices $\xi = \{n_{k}\}_{k=1}^{\infty}$ we introduce a sequence of matrices
\begin{equation}\label{def_A_xi}
B_{1}^{(\xi)} = \vprod\limits_{n=1}^{n_{1}}B_{n},\quad
B_{k+1}^{(\xi)} = \vprod\limits_{n=n_{k}+1}^{n_{k+1}}B_{n},\quad
\forall\, k\in \mathbb{N}.
\end{equation}
Representing $B_{k}^{(\xi)}$ by means of SVD in a form
\begin{equation}\label{svd_A_xi}
B_{k}^{(\xi)} = 
R\Bigl(\varphi_{k}^{(\xi)}\Bigr)Z\Bigl(\lambda_{k}^{(\xi)}\Bigr)R\Bigl(\chi_{k}^{(\xi)}\Bigr)
\end{equation}
we consider a sequence of matrices
\begin{align}\label{def_D_xi}
&D_{k}^{(\xi)} = R\Bigl(\Phi_{k}^{(\xi)}\Bigr)Z\Bigl(\lambda_{k}^{(\xi)}\Bigr),\\
\nonumber
&\Phi_{k}^{(\xi)} = \chi_{k+1}^{(\xi)} + \varphi_{k}^{(\xi)},\quad k\in\mathbb{N}.
\end{align}
Then we arrive at the following
\begin{lems}
Let $\{ B_{n}\}$ be a sequence of $SL(2,\mathbb{R})$-matrices of the form $(\ref{coc_B})$ such that 
\begin{equation*}
\lambda_{n} \ge \Lambda_{0},\quad \forall\, n\in \mathbb{N}
\end{equation*}
with some constant $\Lambda_{0}>1$. Assume that for an increasing sequence of indices 
$\xi = \{n_{k}\}_{k=1}^{\infty}$ a cocycle generated by the sequence $\{D_{k}^{(\xi)}\}$ satisfies conditions of Lemma 3 with constants $\Lambda_{0}^{(\xi)}, C_{\lambda}^{(\xi)}$ and $\delta^{(\xi)}$. If there exist constants $C_{B}>1$ and $N_{0}\in \mathbb{N}$ such that
\begin{align}\label{cond_L4}
\nonumber
&\Vert B_{n}\Vert \le C_{B},\quad \forall\, n\in\mathbb{N},\\
&n_{k+1} - n_{k} \le N_{0},\quad 
\forall\, k\in\mathbb{N},\\
\nonumber
&(\hat C_{\lambda}^{(\xi)})^{\frac{1}{N_{0}}}>1
\end{align}
with
\begin{equation*}
\hat C_{\lambda}^{(\xi)}=
\frac{\Lambda_{0}C_{\lambda}^{(\xi)}}
{\sqrt{\Lambda_{0}^{2}+\left(C_{\lambda}^{(\xi)}\right)^{2}}},
\end{equation*}
then a cocycle generated by the sequence $\{B_{n}\}$ admits estimate $(\ref{coc_est})$.
\end{lems}
{\bf PROOF:} It follows from Lemma 2 that 
\begin{equation}\label{est_m}
m(\lambda_{1}, \lambda_{2}, \varphi)\ge 
\frac{\lambda_{2}}{\lambda_{1}} + \frac{\lambda_{1}}{\lambda_{2}}
\end{equation}
uniformly with respect to $\varphi$.
Note also that 
\begin{equation*}
\left\Vert \vprod\limits_{n=i}^{j}B_{n}\right\Vert \le (C_{B})^{\vert j-i\vert + 1},\quad
\forall\, i,j\in \mathbb{N}.
\end{equation*}
Fix large natural number $k_{0}$. It follows from assumptions and Lemma 3 that for any $k\in \mathbb{N}$ such that $k\ge k_{0}$ one has 
\begin{equation*}
\left\Vert \vprod\limits_{n=1}^{n_{k}+s}B_{n}\right\Vert \ge 
\Lambda_{0}^{(\xi)}\frac{(\hat C_{\lambda}^{(\xi)})^{k}}{(C_{B})^{s}}
\end{equation*}
for any natural $s \le n_{k+1}-n_{k}$. Besides, due to the second condition 
$(\ref{cond_L4})$ one has $n_{k}\le N_{0}k$.
Taking this into account together with the third inequality $(\ref{cond_L4})$, we conclude that
\begin{equation*}
\left\Vert \vprod\limits_{n=1}^{n_{k}+s}B_{n}\right\Vert \ge 
\Lambda_{0}^{(\xi)}\frac{(\hat C_{\lambda}^{(\xi)})^{k}}{(C_{B})^{s}}\ge
\Lambda_{0}^{(\xi)} C_{\Lambda}^{n_{k}+s}, \quad
C_{\Lambda} = \frac{(\hat C_{\lambda}^{(\xi)})^{\frac{k_{0}}{(k_{0}+1)N_{0}}}}
{(C_{B})^{\frac{1}{(k_{0}+1)}}}.
\end{equation*}
Finally, we choose $k_{0}$ to be sufficiently large such that
\begin{equation*}
C_{\Lambda} > 1.
\end{equation*}
It is possible since the third iequality $(\ref{cond_L4})$ is strict. Hence, the cocycle $M_{n}$ generated by the sequence $\{B_{n}\}$ satisfies 
\begin{equation}\label{ineq_M}
\Vert M_{n}\Vert \ge \Lambda_{0}^{(\xi)}\Bigl(C_{\lambda}\Bigr)^{n},\quad
\forall n\ge n_{k_{0}}.
\end{equation}
Decreasing (if necessary) the constant $\Lambda_{0}^{(\xi)}$, one achieves fulfillment of 
$(\ref{ineq_M})$ for all $n\in \mathbb{N}$. This finishes the proof. $\qed$

{\bf Remark} We point out that in the case when a sequence of matrices $\{B_{n}\}$ is generated by a continued fraction $\{(-1, b_{n}), f_{n}\}$, a passage from the sequence 
$\{B_{n}\}$ to a sequence $\{B_{k}^{(\xi)}\}$ via formulae $(\ref{def_A_xi})$ corresponds to a passage from this continued fraction to its contraction 
$\{(\hat a_{k}, \hat b_{k}), \hat f_{k}\}$ (see e.g. \cite{JoTh} for definition) such that
\begin{equation*}
\hat f_{k} = f_{n_{k}},\quad \forall\, k\in\mathbb{N}.
\end{equation*}
The convergence of any contraction of a continued fraction follows from the convergence of the continued fraction. The opposite is not true in general. However, Lemma 4 provides conditions on a contraction, which guarantee the convergence of the continued fraction itself.

The following lemma suggests a receipt how to construct such a contraction. Consider a product $P_{2}$:
\begin{equation*}
P_{2} = R(\varphi_{2})Z(\lambda_{2})R(\varphi_{1})Z(\lambda_{1}).
\end{equation*}
By Lemma 2 it can be represented as
\begin{equation*}
P_{2} = R(\varphi_{2}+\psi)Z(\mu)R(\chi),
\end{equation*}
where $\mu, \psi, \chi$ are defined by $(\ref{f_P1})$.
\begin{lems}
Assume that for positive constants $\Lambda_{0}, C_{\lambda}$ the following conditions hold true
\begin{align}\label{cond_L5_1}
\nonumber
&1.\quad 1 < C_{\lambda} < \Lambda_{0},\\
&2.\quad \lambda_{k}\ge \Lambda_{0}>1, k=1,2,\\
\nonumber
&3.\quad \lambda_{1}>\lambda_{2}.
\end{align}
If, additionally, 
\begin{equation}\label{cond_L5_2}
\vert u\vert \le \varkappa 
\left(\sqrt{\sinh^{2}\left(\ln\frac{\Lambda_{0}}{C_{\lambda}}\right)+\rho} -
\sinh\left(\ln\frac{\Lambda_{0}}{C_{\lambda}}\right)\right),
\end{equation}
where 
\begin{equation}\label{def_u_kappa_rho}
u = \lambda_{2}^{2}\cot \varphi_{1},\quad
\varkappa = \frac{1-\lambda_{1}^{-4}}{1-\lambda_{1}^{-4}\lambda_{2}^{-4}},\quad
\rho = 
\frac{1-(\lambda_{2}^{4}+\lambda_{2}^{-4})\lambda_{1}^{-4}+\lambda_{1}^{-8}}
{1-2\lambda_{1}^{-4}+\lambda_{1}^{-8}}
\end{equation}
and
\begin{equation}\label{cond_L5_3}
\cot\varphi_{2}\in 
\left[\frac{\Lambda_{0} + C_{\lambda}\cot\psi}{\Lambda_{0}\cot\psi - C_{\lambda}},
\frac{\Lambda_{0} - C_{\lambda}\cot\psi}{\Lambda_{0}\cot\psi + C_{\lambda}}\right]
\end{equation}
then
\begin{equation*}
\Lambda_{0}\bigl\vert \cot (\varphi_{2}+\psi)\bigr\vert \ge C_{\lambda}.
\end{equation*} 
\end{lems}
{\bf PROOF:} Similar to Lemma 3 we introduce
\begin{equation*}
\zeta = z(\lambda_{2}, \lambda_{1}, \varphi_{1}).
\end{equation*}
It can be expressed in terms of the variable $u$ as
\begin{equation}\label{f_zeta}
\zeta = \frac{1}{2(1-\lambda_{1}^{-4})}
\left(\Bigl(1-\lambda_{1}^{-4}\lambda_{2}^{-4}\Bigr)u -
\left(1-\frac{\lambda_{2}^{4}}{\lambda_{1}^{4}}\right)u^{-1}\right).
\end{equation} 
Consider an inequality
\begin{equation}\label{ineq_T_psi}
\vert \tan\psi \vert > \frac{\Lambda_{0}}{C_{\lambda}},
\end{equation}
where $\tan\psi$ is defined by $(\ref{f_P1})$.

Due to conditions $(\ref{cond_L5_1})$ 
\begin{equation*}
\vert \tan\psi \vert = \sqrt{1+\zeta^{2}} + \vert\zeta\vert.
\end{equation*}
Thus, a solution of $(\ref{ineq_T_psi})$ in terms of $\vert \zeta\vert$ is 
\begin{equation}\label{ineq_zeta}
\vert \zeta \vert \ge \sinh\left(\ln\frac{\Lambda_{0}}{C_{\lambda}}\right).
\end{equation}
Substituting $(\ref{f_zeta})$ into $(\ref{ineq_zeta})$, one obtains that a solution of 
$(\ref{ineq_T_psi})$ with respect to $\vert u\vert$ is described by $(\ref{cond_L5_2})$.

Finally we note that under condition $(\ref{ineq_T_psi})$ a solution of the inequality
\begin{equation*}
\Lambda_{0}\bigl\vert \cot (\varphi_{2}+\psi)\bigr\vert \ge C_{\lambda}
\end{equation*} 
with respect to $\cot\varphi_{2}$ is given by $(\ref{cond_L5_3})$. $\qed$

{\bf Remark:} In particular, Lemma 5 claims that when both angles $\varphi_{1}$ and $\varphi_{2}$ are close to $\pi/2 \mod \pi$ and, consequently, both multipliers $R(\varphi_{k})Z(\lambda_{k}), k=1, 2$ do not satisfy condition $(\ref{cond_H2_3})$ of Lemma 3, their product $P_{2}$ may satisfy it. 

Indeed, in this case $u\ll 1$ and conditions $(\ref{cond_L5_1})$ imply that 
$\varphi_{1} \approx\frac{\pi}{2}\mod \pi$ and
\begin{equation*}
\zeta = -\frac{\lambda_{1}^{4} - \lambda_{2}^{4}}
{2(\lambda_{1}^{4}-1)} u^{-1} \Bigl(1 + O\bigl(u^{2}\bigr)\Bigr).
\end{equation*}
Taking into account $(\ref{f_P1})$ this leads to
\begin{equation*}
\tan\psi = \frac{\lambda_{1}^{4} - \lambda_{2}^{4}}
{\lambda_{1}^{4}-1} u^{-1} \Bigl(1 + O\bigl(u^{2}\bigr)\Bigr).
\end{equation*}
Hence, inequality $(\ref{ineq_T_psi})$ is fullfilled.

We apply Lemmas 4, 5 to continued fractions and their cocycles generated by a sequence of matrices $\{B_{n}\}$ which do not satisfy hypothesis $(H_{2})$. 

Denote by $\{m_{k}\}$ an increasing sequence of those indices which correspond to violation of $(H_{2})$. Estmates $(\ref{est_m})$ and $(\ref{ineq_M})$ imply that whenever the sequence $\{m_{k}\}$ is finite, the continued fraction still converges. So, in the rest of the paper we will suppose that $\{m_{k}\}$ is infinite. Moreover, we will assume that
\begin{equation}\label{cond_seq_m}
m_{2k} - m_{2k-1} > 1,\quad \forall\, k\ge K_{0}
\end{equation} 
for some $K_{0}\in \mathbb{N}\cup\{0\}$.

Let $\xi$ be a sequence of indices $\{n_{k}\}$ such that
\begin{equation*}
n_{k} = m_{2(k+K_{0})-1}.
\end{equation*}
Applying $(\ref{def_A_xi})$, $(\ref{svd_A_xi})$ and $(\ref{def_D_xi})$ one constructs a sequence of matrices $D_{k}^{(\xi)}$, corresponding to the sequence of indices $\xi$.

Under condition $(\ref{cond_seq_m})$ we consider a sequence of indices $\{l_{k}\}$, which satisfy
\begin{equation}\label{cond_seq_l}
l_{k}\in \left(m_{2(k+K_{0})-1}, m_{2(k+K_{0})}\right),\quad \forall\, k\in \mathbb{N}.
\end{equation} 
Introduce the following sequences of matrices
\begin{equation*}
B_{k}^{(\xi),0} = \vprod\limits_{n=n_{k}+1}^{l_{k}}B_{n},\quad 
B_{k}^{(\xi),-} = \vprod\limits_{n=l_{k}+1}^{j_{k}}B_{n},\quad 
B_{k}^{(\xi),+} = \vprod\limits_{n=j_{k}+1}^{n_{k+1}}B_{n}, 
\end{equation*}
where
\begin{equation*}
j_{k} = m_{2(k+K_{0})},\quad \forall\, k\in\mathbb{N}.
\end{equation*}
By means of SVD we represent
\begin{equation*}
B_{k}^{(\xi), s} = 
R\left(\varphi_{k}^{(\xi),s}\right) 
Z\left(\lambda_{k}^{(\xi),s}\right) R\left(\chi_{k}^{(\xi),s}\right),
\quad s\in \{0, -, +\}
\end{equation*}
and introduce
\begin{align*}
&D_{k}^{(\xi), s} = R\left(\Phi_{k}^{(\xi),s}\right) Z\left(\lambda_{k}^{(\xi),s}\right),\\
\nonumber
&\Phi_{k}^{(\xi),s} = \chi_{k+1}^{(\xi),s} + \varphi_{k}^{(\xi),s},\quad 
\forall\, k\in \mathbb{N}.
\end{align*}
Then one has
\begin{align*}
&B_{k}^{(\xi)} = B_{k}^{(\xi),+}\cdot B_{k}^{(\xi),-}\cdot B_{k}^{(\xi),0},\quad
\forall\, k\in \mathbb{N},\\
&\left\Vert D_{k}^{(\xi)}\right\Vert = 
\left\Vert D_{k}^{(\xi),+}\cdot D_{k}^{(\xi),-}\cdot D_{k}^{(\xi),0}\right\Vert,\quad
\forall\, k\in \mathbb{N}.
\end{align*}
Notice that a product $D_{k}^{(\xi),+}\cdot D_{k}^{(\xi),-}$ has the structure of the product $P_{2}$. Moreover, since $D_{m_{k}}$ does not satisfy condition $(H_{2})$ for all $k\in \mathbb{N}$, one may expect that both of $D_{k}^{(\xi),+},  D_{k}^{(\xi),-}$ also do not satisfy it. We use SVD decomposition to represent
\begin{equation*}
Z\left(\lambda_{k}^{(\xi),+}\right)D_{k}^{(\xi), -} = 
R\left(\Psi_{k}^{(\xi)}\right) Z\left(\mu_{k}^{(\xi)}\right) R\left(\Theta_{k}^{(\xi)}\right)
\end{equation*}
and introduce
\begin{align*}
&u_{k} = \Bigl(\lambda_{k}^{(\xi),+}\Bigr)^{2}\cot\Phi_{k}^{(\xi),-},\quad
\varkappa_{k}=
\frac{1-\Bigl(\lambda_{k}^{(\xi),-}\Bigr)^{-4}}
{1-\Bigl(\lambda_{k}^{(\xi),-}\Bigr)^{-4}\Bigl(\lambda_{k}^{(\xi),+}\Bigr)^{-4}},\\
&\rho_{k} = 
\frac{1 - \left(\Bigl(\lambda_{k}^{(\xi),+}\Bigr)^{4} + 
\Bigl(\lambda_{k}^{(\xi),+}\Bigr)^{-4}\right)\Bigl(\lambda_{k}^{(\xi),-}\Bigr)^{-4} +
\Bigl(\lambda_{k}^{(\xi),-}\Bigr)^{-8}}
{1-2\Bigl(\lambda_{k}^{(\xi),-}\Bigr)^{-4}+\Bigl(\lambda_{k}^{(\xi),-}\Bigr)^{-8}}.
\end{align*}

Then one combines Lemmas 3-5 to formulate the following
\begin{thms} Assume that for a given continued fraction $\{(-1, b_{j}), f_{j}\}$ the condition $(H_{1})$ holds true. Let $\Lambda_{0}^{(\min)}$ be a constant guaranteed by Lemma 1 and $\{B_{n}\}$ be the sequence of matrices (\ref{coc_B}) associated to the continued fraction.
Suppose also that there exist positive constants $C_{\lambda}$ and $\delta$ satisfying (\ref{cond_H2_1}), (\ref{cond_H2_2}) such that 
\begin{equation*}
\Lambda_{0}^{(\min)}\bigl\vert \cot(\Phi_{n})\bigr\vert \ge \delta C_{\lambda},\quad \forall\, n\in\mathbb{N}\setminus \{m_{k}\}_{k=1}^{\infty},
\end{equation*}
where $\{m_{k}\}_{k=1}^{\infty}$ is an increasing sequence of indices.
If there exist positive constants $\Lambda_{0}^{(\max)}>1$, $G_{\lambda}>1$, a natural number $N_{0}$ and a sequence of indices $\{l_{k}\}$ satisfying $(\ref{cond_seq_l})$ such that for all $k\in\mathbb{N}$ the following conditions:
\begin{align}\label{cond_T3}
\nonumber
&1.\quad \lambda_{k}\le \Lambda_{0}^{(\max)},\\
\nonumber
&2.\quad \hat C_{\lambda} = 
\frac{\Lambda_{0}^{(\min)}C_{\lambda}}
{\sqrt{\Bigl(\Lambda_{0}^{(\min)}\Bigr)^{2}+(C_{\lambda})^{2}}} > 1,\\
\nonumber
&3.\quad \hat G_{\lambda} = 
\frac{\Lambda_{0}^{(\min)}G_{\lambda}}
{\sqrt{\Bigl(\Lambda_{0}^{(\min)}\Bigr)^{2}+(G_{\lambda})^{2}}} > 1,\\
&4.\quad n_{k+1} - j_{k}\le N_{0},
\end{align}
%\\
%\nonumber
\begin{align*}
&5.\quad \left(\hat C_{\lambda}\right)^{j_{k}-l_{k}} \ge 
G_{\lambda}\Bigl(\Lambda_{0}^{\max}\Bigr)^{n_{k+1}-j_{k}},\\
\nonumber
&6.\quad \vert u_{k}\vert \le \varkappa_{k} 
\left(\sqrt{\sinh^{2}\left(\ln\frac{\Lambda_{0}^{(\min)}}{G_{\lambda}}\right)+\rho_{k}} -
\sinh\left(\ln\frac{\Lambda_{0}^{(\min)}}{G_{\lambda}}\right)\right),\\
\nonumber
&7.\quad \tan\left(\Phi_{k}^{(\xi),+}\right)\in \left[
\frac{\Lambda_{0} + G_{\lambda} \cot\left(\Psi_{k}^{(\xi)}\right)}
{\Lambda_{0}\cot\left(\Psi_{k}^{(\xi)}\right) - G_{\lambda}},\,\,
\frac{\Lambda_{0} - G_{\lambda} \cot\left(\Psi_{k}^{(\xi)}\right)}
{\Lambda_{0}\cot\left(\Psi_{k}^{(\xi)}\right) + G_{\lambda}}\right]
\end{align*}
hold true, then a cocycle associated with the continued fraction admits estimate 
$(\ref{coc_est})$.
\end{thms}
{\bf PROOF:} First we note that Lemma 3 yields
\begin{eqnarray*}
&\left\Vert D_{k}^{(\xi), 0}\right\Vert \ge 
\Lambda_{0}^{\min}\left(\hat C_{\lambda}\right)^{l_{k}-n_{k}},\quad
\left\Vert D_{k}^{(\xi), -}\right\Vert \ge 
\Lambda_{0}^{\min}\left(\hat C_{\lambda}\right)^{j_{k}-l_{k}},\quad
\left\Vert D_{k}^{(\xi), +}\right\Vert \ge 
\Lambda_{0}^{\min}\left(\hat C_{\lambda}\right)^{n_{k+1}-j_{k}},\\
&\Lambda_{0}\left\vert\cot\left(\Phi_{k}^{(\xi),0}\right)\right\vert \ge C_{\lambda},\quad
\left\Vert D_{k}^{(\xi),+}\right\Vert \le \Bigl(\Lambda_{0}^{\max}\Bigr)^{n_{k+1}-j_{k}}.
\end{eqnarray*}
Hence, due to these estimates together with condtions 3-5 of $(\ref{cond_T3})$ one may apply Lemma 5 to the product $D_{k}^{(\xi),+}\cdot D_{k}^{(\xi),-}$. It gives
\begin{equation*}
\left\Vert D_{k}^{(\xi),+}\cdot D_{k}^{(\xi),-}\right\Vert \ge 
\Lambda_{0}^{\min}\hat G_{\lambda},\quad
\Lambda_{0}^{\min}
\left\vert \cot\left(\Phi_{k}^{(\xi),+}+\Psi_{k}^{(\xi)}\right)\right\vert \ge G_{\lambda}.
\end{equation*}
Then, taking into account the second condition $(\ref{cond_T3})$, we conclude by Lemma 4 that a cocycle generated by matrices $\{B_{n}\}$ with 
$n\in [l_{k}+1, n_{k+1}], k\in \mathbb{N}$ admits estimate $(\ref{coc_est})$. Finally, we notice that inclusion of matrices $\{B_{n}\}$ with $n\in [n_{k}+1, l_{k}], k\in \mathbb{N}$ to the cocycle does not destroy estimate $(\ref{coc_est})$, since all these matrices satisfy conditions of Lemma 3. This finishes the proof. $\qed$

{\bf Remarks:} 1. It has to be noted that a similar theorem is valid if one replaces condition 
$(\ref{cond_seq_m})$ by 
\begin{equation}\label{cond_seq_m_2}
m_{2k+1} - m_{2k} > 1,\quad \forall\, k\ge K_{0}
\end{equation} 
for some $K_{0}\in \mathbb{N}\cup\{0\}$.

2.  We point out that there are two special cases when conditions $(\ref{cond_T3})$ can be simplified. First of them occurs if 
\begin{equation*}
j_{k} - l_{k} \ge K_{0},\,\, \forall\, k\in\mathbb{N},\quad K_{0}\gg 1.
\end{equation*}
In this case the sequence of indices corresponding to violation of condition $(H_{2})$ can be divided into pairs of indices such that the difference between indices in each pair is bounded, but two different pairs are isolated from each other by a large number.

The second case corresponds to
\begin{equation*}
C_{\lambda}\gg 1,\quad \frac{\Lambda_{0}^{\max}}{C_{\lambda}} = o(C_{\lambda}),\quad j_{k}-l_{k}>n_{k+1}-j_{k},\,\, \forall\, k\in\mathbb{N}.
\end{equation*}
In this case the sequence of indices corresponding to violation of condition $(H_{2})$ can be divided into pairs of indices such that the distance between neighbouring pairs can be comparable with the difference between indices constituting a pair, but the cocycle of matrices, which satisfy condition $(H_{2})$, possesses "strong" hyperbolic behaviour. As an example, one may consider a continued fraction $\{(-1, b_{j}), f_{j}\}$ such that
\begin{equation}\label{def_seq_b}
\forall\, j\notin \{t_{k}\}_{k=1}^{\infty}\quad b_{j} = g \hat b_{j},\quad  0 < C_{1} \le \vert\hat b_{j}\vert \le C_{2},\quad g\gg 1,
\end{equation}
where $\{t_{k}\}$ is an increasing sequence of indices and $C_{1,2}$ are constants. Thus, we arrive at the following corollary to Theorem 3.
\begin{cors}
Assume that for a given increasing sequence of natural numbers $\{t_{k}\}$ a continued fraction $\{(-1, b_{j}), f_{j}\}$ satisfies conditions $(\ref{def_seq_b})$, $(H_{1})$ and 
\begin{align}\label{cond_C2_1}
\nonumber
&1.\quad \exists\, K_{0}\in\mathbb{N}\cup\{0\}:\quad
t_{k+1}-t_{k} \ge 4,\quad \forall\, k\ge K_{0} ,\\
&2.\quad \left\vert b_{j} - \frac{b_{j-1}}{1+b_{j-1}^{2}}\right\vert < 
\frac{1}{2C_{1}^{2}g^{2}},\quad \forall\, j\in \{t_{k}\}:\,\, j = 2s,\, s\in\mathbb{N}\\
\nonumber
&3.\quad \left\vert b_{j} - \frac{b_{j+1}}{1+b_{j+1}^{2}}\right\vert < 
\frac{1}{2C_{1}^{2}g^{2}},\quad \forall\, j\in \{t_{k}\}:\,\, j = 2s-1,\, s\in\mathbb{N}.
\end{align}
If for a sequence $\xi=\{m_{k}\}$ such that
\begin{equation}\label{def_seq_m_2}
m_{k} = \left[ \frac{t_{k}}{2}\right],\quad \forall\, k\in \mathbb{N},
\end{equation}
where $[x]$ stands for the integer part of a real number $x$,
there exist a natural number $N_{0}$, positive constant $\alpha$ and a sequence of indices $\{l_{k}\}$ satisfying 
$(\ref{cond_seq_l})$ such that for all $k\in\mathbb{N}$ the following conditions:
\begin{align}\label{cond_C2_2}
\nonumber
&4.\quad n_{k+1}-j_{k}\le N_{0},\\
&5.\quad j_{k}-l_{k} > n_{k+1}-j_{k},\\
\nonumber
&6.\quad \left\vert \cot \Phi_{k}^{(\xi),-}\right\vert \le 
\Bigl(C_{2}g\Bigr)^{-4(n_{k+1}-j_{k})+2-\alpha},\\ 
\nonumber
&7.\quad \cot \Phi_{k}^{(\xi),+}\cdot \cot \Phi_{k}^{(\xi),-} < 0,
\end{align}
hold true, then there exists $g_{0}$ such that for all $g>g_{0}$ the continued fraction $\{(-1, b_{j}), f_{j}\}$ converges.
\end{cors}

{\bf PROOF:} It has to be noted that the first condition $(\ref{cond_C2_1})$ yields both conditions $(\ref{cond_seq_m})$, $(\ref{cond_seq_m_2})$ for the sequence 
$(\ref{def_seq_m_2})$.

Consider formulae $(\ref{SVD_2})$ and assume that $2n-1, 2n\notin \{t_{j}\}$. In this case
\begin{equation*}
\vert b_{2n-1}\vert \ge C_{1}g,\quad \vert b_{2n}\vert \ge C_{1}g,
\end{equation*}
which leads to
\begin{equation}\label{est_ctg_phi_0}
\lambda_{n}\ge \Bigl(C_{1}g\Bigr)^{2},\quad
\vert \cot\Phi_{n}\vert \ge C_{1}g (1+O(g^{-2})).
\end{equation}
Thus, matrices $B_{n}$ satisfy conditions of Lemma 3, for example, with constants 
\begin{equation}\label{def_const}
\Lambda_{0}^{(\min)} = C_{1}g,\quad 
C_{\lambda}=\Lambda_{0}^{(\min)}/2,\quad \delta = 2
\end{equation} 
provided $g$ to be sufficiently large.

Assume now that $2n\in \{t_{j}\}$. As it was mentioned in Lemma 1, equality 
$(\ref{eq_b_2k})$ implies
\begin{equation*}
\cot\varphi_{n} = 0. 
\end{equation*}
Expanding $\cot\varphi_{n}$ near 
\begin{equation*}
b_{2n}^{*} = \frac{b_{2n-1}}{1+b_{2n-1}^{2}},
\end{equation*}
one obtains
\begin{equation*}\label{Taylor_ctg_phi}
\vert \cot\varphi_{n}\vert = 
\frac{\left(1+b_{2n-1}^{2}\right)^{2}}
{b_{2n-1}^{2} \left(b_{2n-1}^{2}+2\right)}\cdot
\Bigl\vert b_{2n} - b_{2n}^{*}\Bigr\vert\cdot
\left(1 + O\left(\Bigl\vert b_{2n} - b_{2n}^{*}\Bigr\vert\right)\right).
\end{equation*}
Since $2n-1\notin \{t_{j}\}$ then $\vert b_{2n-1}\vert \ge C_{1}g$ and we get
\begin{equation}\label{est_ctg_phi_1}
\lambda_{n}\ge C_{1}g \left(1 + O(g^{-2})\right),\quad
\vert \cot\Phi_{n}\vert = C_{1}g\cdot \Bigl\vert b_{2n} - b_{2n}^{*}\Bigr\vert\cdot
\left(1 + O\left(g^{-2} + \Bigl\vert b_{2n} - b_{2n}^{*}\Bigr\vert\right)\right).
\end{equation}
Thus, the second condition $(\ref{cond_C2_1})$ guarantees violation of $(H_{2})$ for $B_{n}$ with $2n\in \{t_{j}\}$.

The case $2n-1\in \{t_{j}\}$ can be considered in a similar way. Then one gets 
\begin{equation}\label{est_ctg_phi_2}
\lambda_{n-1}\ge C_{1}g \left(1 + O(g^{-2})\right),\quad
\vert \cot\Phi_{n-1}\vert = 
C_{1}g\cdot \Bigl\vert b_{2n-1} - b_{2n-1}^{*}\Bigr\vert\cdot
\left(1 + O\left(g^{-2} + \Bigl\vert b_{2n-1} - b_{2n-1}^{*}\Bigr\vert\right)\right),
\end{equation}
where
\begin{equation*}
b_{2n-1}^{*} = \frac{b_{2n}}{1+b_{2n}^{2}}.
\end{equation*}
We conclude that the third condition $(\ref{cond_C2_1})$ leads to violation of $(H_{2})$ for $B_{n-1}$ with $2n-1\in \{t_{j}\}$. 

It follows from $(\ref{SVD_2})$ that whenever $2n-1, 2n \notin \{t_{j}\}$
\begin{equation}\label{est_lambda_0}
\lambda_{n}\le (C_{2}g)^{2}(1+O(g^{-2})).
\end{equation}
If $2n\in \{t_{j}\}$ we expand $\lambda_{n}^{2}$ near $b_{2n}=b_{2n}^{*}$:
\begin{equation*}
\lambda_{n}^{2} = 1 + b_{2n-1}^{2} + 
\frac{1}{2}\left(1+b_{2n-1}^{2}\right)^{2}
\left(\frac{1}{b_{2n-1}^{2}+2}+\frac{1}{b_{2n-1}^{2}}\right)
\Bigl( b_{2n} - b_{2n}^{*}\Bigr)^{2}\cdot
\left(1 + O\left(\Bigl\vert b_{2n} - b_{2n}^{*}\Bigr\vert\right)\right)
\end{equation*}
Then the second condition $(\ref{cond_C2_1})$ gives
\begin{equation}\label{est_lambda_1}
\lambda_{n}\le C_{2}g(1+O(g^{-2})).
\end{equation}
In a similar way one obtains that for $2n-1\in \{t_{j}\}$ 
\begin{equation}\label{est_lambda_2}
\lambda_{n-1}\le C_{2}g(1+O(g^{-2})).
\end{equation}
We combine estimates $(\ref{est_ctg_phi_0})$, $(\ref{est_ctg_phi_1})$, 
$(\ref{est_ctg_phi_2})$ together with $(\ref{est_lambda_0})$, $(\ref{est_lambda_1})$, 
$(\ref{est_lambda_2})$ and get
\begin{align}\label{est_lambda_xi}
&\left(\hat C_{1} g\right)^{2(j_{k}-l_{k})-1}\le \lambda_{k}^{(\xi),-} \le
\left(\hat C_{2} g\right)^{2(j_{k}-l_{k})-1},\\
\nonumber
&\left(\hat C_{1} g\right)^{2(n_{k+1}-j_{k})-1}\le \lambda_{k}^{(\xi),+} \le
\left(\hat C_{2} g\right)^{2(n_{k+1}-j_{k})-1},
\end{align}
where $\hat C_{1} < C_{1}$, $\hat C_{2} > C_{2}$ are positive constants. Note that 
$\hat C_{1}, \hat C_{2}$ can be chosen arbitrary close to $C_{1}, C_{2}$, respectively, for sufficiently large $g$.
Hence
\begin{equation}\label{est_frac_lambda}
\frac{\lambda_{k}^{(\xi),-}}{\lambda_{k}^{(\xi),+}} \ge
\frac{\left(\hat C_{1}\right)^{2(j_{k}-l_{k})-1}}
{\left(\hat C_{2}\right)^{2(n_{k+1}-j_{k})-1}}
g^{2(2j_{k}-n_{k+1}-l_{k})} = 
\frac{\left(\hat C_{1}g^{(1-\eta)}\right)^{2(j_{k}-l_{k})}}
{\hat C_{1}\left(\hat C_{2}\right)^{2(n_{k+1}-j_{k})-1}}
g^{2\bigl(\eta(j_{k}-l_{k})-(n_{k+1}-j_{k})\bigr)}
\end{equation}
for any $0 < \eta < 1$.
We set 
\begin{equation*}
\eta = \frac{N_{0}+\nu}{N_{0} + 1},\quad 0<\nu<1.
\end{equation*}
Then, due to $(\ref{est_frac_lambda})$ and the second inequality $(\ref{cond_C2_2})$ one obtains
\begin{equation*}
\frac{\lambda_{k}^{(\xi),-}}{\lambda_{k}^{(\xi),+}} \ge
\frac{\left(\hat C_{1}g^{(1-\eta)}\right)^{2(j_{k}-l_{k})}}
{\hat C_{1}\left(\hat C_{2}\right)^{2(n_{k+1}-j_{k})-1}}
g^{2\nu} \ge
g^{2\nu}
\end{equation*}
provided $g$ to be sufficiently large. This allows to define
\begin{equation*}
G_{\lambda} = g^{2\nu},\quad 0<\nu<1.
\end{equation*}
We consider $\nu$ such that $\nu<1/2$. In this case
\begin{equation}\label{def_G_lambda}
\frac{\Lambda_{0}^{(\min)}}{G_{\lambda}} = C_{1}g^{1-2\nu}\gg 1,
\end{equation}
hence, it is possible to apply Lemma 5.

Besides, one has
\begin{equation}\label{est_kappa_rho}
\varkappa_{k} = 1 + O(g^{-8}),\quad \rho_{k} = 1 + O(g^{-8})
\end{equation}

Substituting $(\ref{est_lambda_xi})$, $(\ref{def_G_lambda})$, $(\ref{est_kappa_rho})$ into inequality 6 of $(\ref{cond_T3})$, we obtain that an estimate
\begin{equation*}
\left(\hat C_{2} g\right)^{4(n_{k+1}-j_{k})-2}\Bigl\vert \cot \Phi_{k}^{(\xi),-}\Bigr\vert \le
\frac{1}{2} C_{1}^{-1}g^{-(1-2\nu)} = \frac{G_{\lambda}}{2\Lambda_{0}}
\end{equation*}
implies condition 6 in $(\ref{cond_T3})$. But this estimate follows from condition 6 in 
$(\ref{cond_C2_2})$ if one fixes $\nu$ such that
\begin{equation*}
\nu > \frac{1-\alpha}{2}.
\end{equation*}
Finally, we notice that
\begin{equation*}
\Phi_{k}^{(\xi),+} = \Phi_{m_{2k+1}} + O(g^{-4}),\quad
\left\vert\cot \Psi_{k}^{(\xi)}\right\vert \le 
C_{2}^{-1}g^{-\alpha}\left(1+O\left(g^{-\alpha}\right)\right) \le 
\frac{G_{\lambda}}{\Lambda_{0}}.
\end{equation*}
Taking this into account together with $(\ref{def_G_lambda})$,  one obtains that condition 7 in $(\ref{cond_T3})$ is satisfied provided condition 7 of $(\ref{cond_C2_2})$ holds true. This finishes the proof. $\qed$

{\bf Remark:} It has to be noted that although 
\begin{equation*}
\Phi_{k}^{(\xi),-} = \Phi_{m_{2k}} + O(g^{-4})
\end{equation*}
one cannot replace $\Phi_{k}^{(\xi),-}$ by $\Phi_{m_{2k}}$ in condition 6 of 
$(\ref{cond_C2_2})$.

\section{Critical set method}

In the rest of the paper we study convergence of a functional continued fractions generated by a pair $(b, \omega)$ with $b\in C(\mathbb{T}^{1}, \mathbb{R})$, 
$\omega\in (0, 1)\setminus \mathbb{Q}$ in the sense of Definition 2. We suppose that a function $b$ satisfies condition $(H'_{1})$, but condition $(H'_{2})$ is violated.

More precisely, we consider a functional continued fraction, generated by a sequence
\begin{equation*}\label{def_b(x)}
b_{j}(x) = g b(x-(j-1)\omega),\quad j\in \mathbb{Z}
\end{equation*}
with $b\in C^{1}(\mathbb{T}^{1},\mathbb{R})$, $\omega\in (0, 1)\setminus \mathbb{Q}$ and $g\gg 1$.

It is supposed the following condition is fullfilled
\begin{equation*}
(H''_{1})\quad\quad
\min\limits_{x\in \mathbb{T}^{1}} \biggl( (b(x)-b(x-\omega))^{2}+
g^{2}b^{2}(x) b^{2}(x-\omega)\biggr) \ge h_{*}^{2} > 0,
\end{equation*}
where $h_{*}$ does not depend on $g$.

The associated cocycle $(\ref{coc_B_2})$ can be constructed by means of linear transformation
\begin{align}\label{coc_B_3}
&B(x) = R(\Phi(x))Z(\lambda(x)),\\
\nonumber
&\Phi(x) = \psi(x-\omega) + \varphi(x),
\end{align}
where expressions for the functions $\lambda, \varphi, \chi$ are given by formulae 
$(\ref{SVD_2})$ with $b_{2k}, b_{2k-1}$ being replaced by $gb(x-\omega), gb(x)$, respectively.

Following \cite{Iva21} we introduce the critical set 
\begin{equation*}
\mathcal{C}_{0} = \left\{x\in \mathbb{T}^{1}: \cos\Phi(x) = 0\right\}
\end{equation*}
and assume that
\begin{equation*}
(H_{3})\quad\quad\quad 
\mathcal{C}_{0} = \bigcup\limits_{p=0}^{P}\{x_{p}\},\quad
P\in \mathbb{N};\quad
\Phi'(x_{p})\neq 0,\,\, \forall\, p\in \{0,\ldots P\}.
\end{equation*}
It has to be noted that due to condition $(H''_{1})$ and formulae $(\ref{SVD_2})$ 
$\cot \varphi(x)$ and $\cot \psi(x-\omega)$ cannot vanish simultaneously. Moreover, 
$(H''_{1})$ imples that there exists a positive constant $\Lambda_{0}>1$ such that
\begin{equation*}
\lambda(x)\ge \Lambda_{0}g,\quad \forall\, x\in \mathbb{T}^{1}.
\end{equation*}    
Hence, if $\cot \varphi(x) = o(1)$ (resp. $\cot \psi(x-\omega)) = o(1)$) then $\vert \cot \psi(x-\omega)\vert \ge C g$ (resp. $\vert\cot \varphi(x)\vert \ge C g$) with some positive constant $C$. Due to differentiabilty of the functions $\varphi$ and $\psi$, it means that each point of the critical set belongs to a $g^{-1}$-neighbourhood of a solution of the following equation
\begin{equation}\label{eq_Crit_1}
\cot \varphi(x)\cdot \cot \psi(x-\omega) = 0.
\end{equation}
Note that by Lemma 1 equation $(\ref{eq_Crit_1})$ can be rewritten as
\begin{equation}\label{eq_Crit_1_a}
\left(b(x-\omega) - \frac{b(x)}{1+g^{2}b^{2}(x)}\right) \cdot
\left(b(x-\omega) - \frac{b(x-2\omega)}{1+g^{2}b^{2}(x-2\omega)}\right) = 0.
\end{equation}
On the other hand, by the implicit function theorem, each solution of $(\ref{eq_Crit_1_a})$ lies in a $g^{-2}$-neighborhood of a point $x_{p}^{(0)}$, which is a solution of the equation
\begin{equation*}
b(x-\omega) = 0.
\end{equation*}
Thus, condition $(H_{3})$ can be replaced by
\begin{equation*}
(H'_{3})\quad\quad\quad 
\mathcal{C}_{0}^{(0)}  = \bigcup\limits_{p=0}^{P}\{x_{p}^{(0)}\},\quad
P\in \mathbb{N};\quad
b'(x_{p}^{(0)})\neq 0,\,\, \forall\, p\in \{0,\ldots P\},
\end{equation*}
where $\mathcal{C}_{0}^{(0)} = \left\{x\in \mathbb{T}^{1}: b(x-\omega)=0\right\}$.

One notes that condition $(H'_{3})$ implies that $P$ is odd, i.e. the critical set consists of even number of points. We also emphasize that differentiability of the function $b$ can be replaced by its differentiability in small neighborhoods of the critical set $\mathcal{C}_{0}$, its image and preimage under rotation over an angle $-2\pi\omega$. However, to simplify exposition we consider the case $b\in C^{1}(\mathbb{T}^{1},\mathbb{R})$. 

Besides, to shorten formulae in the rest of the paper we introduce for any $\varpi\in \mathbb{R}$ the rotation map over an angle $-2\pi \varpi$ by
\begin{equation*}
\sigma_{\varpi}:\mathbb{T}^{1}\to \mathbb{T}^{1},\quad
\sigma_{\varpi}: x\mapsto x - \varpi \mod 1.
\end{equation*} 
Hence, 
\begin{equation*}
\sigma_{\varpi}^{k}: x\mapsto x - \varpi k \mod 1,\quad
\forall\, k\in\mathbb{Z}.
\end{equation*}

Define for any $p\in \{0,\ldots P\}$ and $\delta>0$ the time of collision 
(see \cite{Iva21})
\begin{equation*}\label{def_T_time}
T_{\delta}(x_{p})=\min\bigl\{k\in \mathbb{N}: 
U_{\delta}\bigl(\sigma^{k}_{2\omega}(x_{p})\bigr) \cap 
U_{\delta}(\mathcal{C}_{0}) \neq \emptyset \bigr\},
\end{equation*}
where $U_{\delta}(x)$ and $U_{\delta}(A)$ stand for $\delta$-neighborhoods of 
$x\in \mathbb{T}^{1}$ and $A\subset \mathbb{T}^{1}$, respectively.

Taking $\delta$ to be sufficiently small, we may assume that there exists a unique index 
$R_{\delta}(p)$ such that 
\begin{equation*}\label{def_r(p)}
\sigma_{2\omega}^{T_{\delta}(x_{p})}(x_{p}) \in U_{\delta}(x_{R_{\delta}(p)}).
\end{equation*}
We will call the case $R_{\delta}(p) = p$ the primary collision, while the case 
$R_{\delta}(p)\neq p$ will be called the secondary collision (see \cite{Iva23}). It has to be noted that the time of primary collision depends only on arithmetical properties of the rotation number $\omega$. In contrary, the time of secondary collision depends also on the function $b$ and, in particular, on distances 
$dist(x_{p}, x_{p'})$, $p\neq p'$, where 
\begin{equation*}
dist(x,y) = \{\vert x - y\vert\},\quad
\forall\, x,y\in \mathbb{T}^{1}
\end{equation*}
and $\{x\}$ denotes the fractional part of $x$.

For each $p = 0,\ldots, P$ we define a period
\begin{equation*}\label{def_S}
S(p) = \min\{k\in \mathbb{N}:\,\, R_{\delta}^{k}(p) = p\}.
\end{equation*}
\begin{defs}
We say that critical points $x_{p}, x_{p'}, p\neq p'$ are cyclically connected if there exists an integer $k < S(p)$ such that
\begin{equation*}
R_{\delta}^{k}(p) = p'.
\end{equation*}
\end{defs}
Notice that if $x_{p}, x_{p'}, p\neq p'$ are cyclically connected then $S(p) = S(p')$. Besides, we remark that cyclic connectness is an equivalence relation and we denote the class of equivalence for a point $x_{p}\in \mathcal{C}_{0}$ by $[x_{p}]$. Thus, the critical set is divided into classes of equivalence. Let $M$ be the number of such equivalence class. Then enumerating them, one obtains
\begin{equation*}
\mathcal{C}_{0} = \bigcup\limits_{i=1}^{M} L_{i},
\end{equation*}
where $L_{i}$ stands for the $i$-th class of equivalence. Since all points inside a class 
$L_{i}$ have the same period, we denote this period by $S_{i}$.

In each class $L_{i}$ we fix a point $x_{p_{*}}$ and denote it by $x_{i, 1}$. Then other points from this class of equivalence can be enumerated such that
\begin{equation*}
\sigma_{2\omega}^{T_{\delta}(x_{i,k})}(x_{i,k})\in U_{\delta}(x_{i,k+1}),\quad
k=1,\ldots, S_{i}-2.
\end{equation*}
This yields
\begin{equation*}
R_{\delta}\bigl((i,k)\bigr) = (i, k+1 \mod S_{i}).
\end{equation*}
In the rest we will omit the notation $"\mod S_{i}"$, but will always presume that the second index in a pair $(i,k)$ is calculated modulo $S_{i}$.
 
One notes that due to conditions $(H''_{1})$ and $(H'_{3})$ the following estimates are valid for sufficiently small $\delta$ 
\begin{align}\label{cond_b}
\nonumber
&\vert b(x) \vert \ge C_{1}\delta,\quad
\forall\, x\in \mathbb{T}^{1}\setminus U_{\delta}(\mathcal{C}_{0}),\\
&b(x) = b'(x_{i,k}) (x-x_{i,k})\cdot (1 + o(1)),\quad
\forall\, x\in U_{\delta}(\mathcal{C}_{0}),\\
\nonumber
&b'(x) = b'(x_{i,k})\cdot (1 + o(1)),\quad
\forall\, x\in U_{\delta}(\mathcal{C}_{0}),
\end{align}
where constant $C_{1}$ can be chosen arbitrary close to 
$\min\limits_{(i,k)}\bigl\vert b'(x_{i,k}^{(0)})\bigr\vert$ by decreasing of $\delta$.

We assign to each class of equivalence $L_{i}$ a set of matrices 
\begin{equation*}
D_{i, k}(x) = \vprod\limits_{j=1}^{T_{\delta}(x_{i,k})}
B\left(\sigma_{2\omega}^{j}(x)\right),\quad
k = 1,\ldots, S_{i}-1,\quad x\in U_{\delta}(x_{i,k}),
\end{equation*}
and represent them by means of SVD decomposition
\begin{equation*}
D_{i, k}(x) = R(\varphi_{i,k}(x)) Z(\lambda_{i,k}(x) R(\chi_{i,k}(x)),\quad
k = 1,\ldots, S_{i}-1,\quad x\in U_{\delta}(x_{i,k}).
\end{equation*}
Define
\begin{equation*}
\Phi_{i,k}(x) = \chi_{i,k+1}\left(\sigma_{2\omega}^{T_{\delta}(x_{i,k})}(x)\right) + \varphi_{i,k}(x),
\quad x\in U_{\delta}(x_{i,k}).
\end{equation*}
The following lemma holds true.
\begin{lems}
There exists a unique solution $x_{i,k}^{(*)}$ of an equation
\begin{equation}\label{eq_Phi_ik}
\cos \Phi_{i,k}(x) = 0,\quad x\in U_{\delta}(x_{i,k}).
\end{equation}
It satisfies an estimate
\begin{equation}\label{est_x_ik}
dist\left(\sigma_{2\omega}^{T_{\delta}(x_{i,k})}\left(x_{i,k}^{(*)}\right), x_{i,k+1}\right) = O\bigl(g^{-2}\bigr).
\end{equation}
\end{lems}
{\bf PROOF:} First we notice that 
\begin{equation*}
\sigma_{2\omega}^{j}(x)\notin U_{\delta}(\mathcal{C}_{0}),\quad
\forall\, j=1,\ldots, T_{\delta}(x_{i,k})-1,\quad 
\forall\, x\in U_{\delta}(x_{i,k}).
\end{equation*}
This leads to an estimate
\begin{equation}\label{est_lambda_good}
\lambda\left(\sigma_{2\omega}^{j}(x)\right)\ge (C_{1}\delta g)^{2},\quad
\forall\, j=1,\ldots, T_{\delta}(x_{i,k})-1,\quad 
\forall\, x\in U_{\delta}(x_{i,k}).
\end{equation}
On the other hand, if $j=T_{\delta}(x_{i,k})$ then it may happen that for some
$x\in U_{\delta}(x_{i,k})$ its image under rotation over the angle $-2\pi T_{\delta}(x_{i,k})$ falls into $U_{\delta}(x_{i,k+1})$ and estimate $(\ref{est_lambda_good})$ is violated. However, due to condition $(H'_{1})$ one has
\begin{equation*}
\lambda\left(\sigma_{2\omega}^{T_{\delta}(x_{i,k})}(x)\right) \ge \Lambda_{0}g,\quad 
\forall\, x\in U_{\delta}(x_{i,k}).
\end{equation*}
Applying consequitively $(\ref{SVD_2})$ to products of matrices 
\begin{equation*}
\vprod\limits_{j=1}^{s}B\left(\sigma_{2\omega}^{j}(x)\right),\quad 
x\in U_{\delta}(x_{i,k})
\end{equation*}
starting from $s=2$ to $s=T_{\delta}(x_{i,k})$, we obtain that
\begin{align}\label{f_Phi_ik}
\nonumber
&\varphi_{i,k}(x) = \varphi\left(\sigma_{2\omega}^{T_{\delta}(x_{i,k})}(x)\right) + \psi_{i,k}(x),\quad
\psi_{i,k}(x) = \sum\limits_{j=1}^{T_{\delta}(x_{i,k})}\psi_{i,k}^{(j)}(x),\quad
\chi_{i,k}(x) = \sum\limits_{j=1}^{T_{\delta}(x_{i,k})}\chi_{i,k}^{(j)}(x),\\
&\psi_{i,k}^{(1)}(x) = O\bigl(g^{-2}\bigr),\quad 
\psi_{i,k}^{(j)}(x) = O\left(g^{-2(2j-1)}\right),\quad
\forall\, j=2,\ldots, T_{\delta}(x_{i,k}),\quad
\forall\, x\in U_{\delta}(x_{i,k}),\\
\nonumber
&\chi_{i,k}^{(1)}(x) = O\bigl(g^{-2}\bigr),\quad 
\chi_{i,k}^{(j)}(x) = O\left(g^{-2(2j-1)}\right),\quad
\forall\, j=2,\ldots, T_{\delta}(x_{i,k}),\quad
\forall\, x\in U_{\delta}(x_{i,k}).
\end{align}
Moreover, it follows from differentiability of the function $b$ and formulae $(\ref{SVD_2})$ that
\begin{equation}\label{est_der_Phi_ik}
\psi'_{i,k}(x) = O\bigl(g^{-2}\bigr),\quad
\chi'_{i,k}(x) = O\bigl(g^{-2}\bigr),\quad
\forall\, x\in U_{\delta}(x_{i,k}).
\end{equation}
Taking into account $(\ref{cond_b})$, $(\ref{f_Phi_ik})$, $(\ref{est_der_Phi_ik})$, we conclude by the implicit function theorem that for sufficiently large $g$ there exists a unique solution of equation $(\ref{eq_Phi_ik})$ which satisfies estimate $(\ref{est_x_ik})$.
$\qed$

As a consequence of Lemma 6 one obtains that there exist positive constants 
$r_{i,k}^{\min}, r_{i,k}^{\max}$ such that
\begin{equation}\label{est_r_ik}
r_{i,k}^{\min} g \left\vert x-x_{i,k}^{(*)}\right\vert \le 
\left\vert\cot \Phi_{i,k}(x)\right\vert \le 
r_{i,k}^{\max} g \left\vert x-x_{i,k}^{(*)}\right\vert. 
\end{equation}
Note that constants $r_{i,k}^{\min}, r_{i,k}^{\max}$ can be chosen close to 
$\vert b'(x_{i,k}^{(0)})\vert$ provided $\delta$ to be sufficiently small.

For each $(i,k), i=1,\ldots, M, k=1,\ldots, S_{i}-1$ we define
\begin{equation*}
\Delta_{i,k} = \sigma_{2\omega}^{T_{\delta}(x_{i,k})}\bigl(x_{i,k}^{(*)}\bigr) - 
x_{i,k+1}^{(*)}.
\end{equation*}
Consider a product
\begin{equation*} 
P_{i,k}^{(2)}(x) = 
R\Bigl(\Phi_{i,k+1}\bigl(\sigma_{2\omega}^{-T_{\delta}(x_{i,k})}(x) \bigr)\Bigr)\cdot Z\Bigl(\lambda_{i,k+1}\bigl(\sigma_{2\omega}^{-T_{\delta}(x_{i,k})}(x)\bigr)\Bigr)\cdot
R\bigl(\Phi_{i,k}(x)\bigr)\cdot Z\bigl(\Phi_{i,k}(x)\bigr)
\end{equation*}
for
\begin{equation*}
x\in U_{\delta}(x_{i,k})\cap 
\sigma_{2\omega}^{-T_{\delta}(x_{i,k})}\bigl(U_{\delta}(x_{i,k+1})\bigr)
\end{equation*}
and represent it in a form
\begin{equation*}
P_{i,k}^{(2)}(x) = 
R\bigl(\varphi_{i,k}^{(2)}(x)\bigr) 
Z\bigl(\lambda_{i,k}^{(2)}(x)\bigr) 
R\bigl(\chi_{i,k}^{(2)}(x)\bigr),\quad 
x\in U_{\delta}(x_{i,k})\cap 
\sigma_{2\omega}^{-T_{\delta}(x_{i,k})}\bigl(U_{\delta}(x_{i,k+1})\bigr).
\end{equation*}
Then we arrive at
\begin{lems}
Suppose that for a given $(i,k), i\in \{1,\ldots, M\}, k\in \{1,\ldots, S_{i}-1\}$ the following conditions are satisfied
\begin{equation}\label{cond_L7}
T_{\delta}(x_{i,k}) > T_{\delta}(x_{i,k+1}),\quad
b'\left(x_{i,k}^{(0)}\right)\cdot b'\left(x_{i,k+1}^{(0)}\right) < 0,\quad
\left\vert \Delta_{i,j}\right\vert \le 
\frac{\left(r_{i,k}^{\max} \cdot r_{i,k+1}^{\max}\right)^{-1/2}}
{g(C_{2}g)^{2 T_{\delta}(x_{i,k+1})-1}},
\end{equation}
where
\begin{equation*}
C_{2}=\max\limits_{x\in \mathbb{T}^{1}}\vert b(x)\vert.
\end{equation*}
Then 
\begin{equation*}
\left\vert \cot\Phi_{i,k}^{(2)}(x)\right\vert \ge 
\frac{\left(r_{i,k+1}^{\max}/r_{i,k}^{\max}\right)^{1/2}}
{(C_{2}g)^{2T_{\delta}(x_{i,k+1})-1}}\cdot
\left(1 + O\left(g^{-4T_{\delta}(x_{i,k+1})+2} + 
g^{-8\bigl(T_{\delta}(x_{i,k})-T_{\delta}(x_{i,k+1})\bigr)}\right)\right).
\end{equation*}
\end{lems}
{\bf PROOF:} Note that estimate $(\ref{est_lambda_good})$ and condition $(H'_{1})$ imply
\begin{equation}\label{est_lambda_ik}
\Lambda_{0}g\left(C_{1}\delta g\right)^{2(T_{\delta}(x_{i,k})-1)} \le \lambda_{i,k}(x) \le
\left(C_{2} g\right)^{2T_{\delta}(x_{i,k})-1},\quad
x\in U_{\delta}(x_{i,k})\cap 
\sigma_{2\omega}^{-T_{\delta}(x_{i,k})}\bigl(U_{\delta}(x_{i,k+1})\bigr).
\end{equation}
Hence, the first condition $(\ref{cond_L7})$ guarantees that
\begin{equation*}
\lambda_{i,k}(x) \gg \lambda_{i,k+1}(x)
\end{equation*}
uniformly with respect to 
$x\in U_{\delta}(x_{i,k})\cap 
\sigma_{2\omega}^{-T_{\delta}(x_{i,k})}\bigl(U_{\delta}(x_{i,k+1})\bigr)$.

Due to estimates $(\ref{est_r_ik})$ the rest of the proof can be done similar to Lemma 4 from \cite{Iva23}. We outline the main idea.

Introduce a function
\begin{equation*}
\zeta_{i,k}(x) = 
z\Bigl(\lambda_{i,k+1}\bigl(\sigma_{2\omega T_{\delta}(x_{i,k})}^{-1}(x)\bigr), \lambda_{i,k}(x), \varphi_{i,k}(x)\Bigr),
\end{equation*}
where the function $z$ is defined by $(\ref{f_P1})$. Then, by $(\ref{f_zeta})$ and 
$(\ref{est_lambda_ik})$ one obtains a uniform estimate
\begin{equation*}
\cot \psi_{i,k}(x) = u_{i,k}(x) \left(1 + 
O\left(g^{-8\bigl(T_{\delta}(x_{i,k})-T_{\delta}(x_{i,k+1})\bigr)}\right)\right),\quad
x\in U_{\delta}(x_{i,k})\cap 
\sigma_{2\omega}^{-T_{\delta}(x_{i,k})}\bigl(U_{\delta}(x_{i,k+1})\bigr).
\end{equation*}
where $\psi_{i,k}(x) = \varphi_{i,k}^{(2)}(x) - 
\Phi_{i,k+1}\bigl(\sigma_{2\omega}^{-T_{\delta}(x_{i,k})}(x)\bigr)$ and 
$u_{i,k}(x) = \lambda_{i,k+1}^{2}\bigl(\sigma_{2\omega}^{-T_{\delta}(x_{i,k})}(x)\bigr)
\cot \Phi_{i,k}(x)$.

We apply $(\ref{est_r_ik})$ to get a similar estimate for $\cot \psi_{i,k}(x)$:
\begin{equation*}\label{est_cot_psi_ik}
r_{i,k}^{\min} g \mu_{i,k}^{2}(x)\left\vert x-x_{i,k}^{(*)}\right\vert \le 
\left\vert\cot \psi_{i,k}(x)\right\vert \le 
r_{i,k}^{\max} g \mu_{i,k}^{2}(x)\left\vert x-x_{i,k}^{(*)}\right\vert
\end{equation*}
with $\mu_{i,k}(x) = 
\lambda_{i,k+1}\bigl(\sigma_{2\omega}^{-T_{\delta}(x_{i,k})}(x)\bigr)$.

Thus, both $\cot \psi_{i,k}(x)$ and 
$\cot \Phi_{i,k+1}\bigl(\sigma_{2\omega}^{-T_{\delta}(x_{i,k})}(x)\bigr)$ are well approximated (in the sense of $(\ref{est_cot_psi_ik})$ and $(\ref{est_r_ik})$, respectively) near $x_{i,k}^{*}$ by linear functions. It leads to
\begin{equation}\label{est_Phi_ik_2}
\left\vert \cot\Phi_{i,k}^{(2)}(x)\right\vert \ge \frac{ga}{1-a\mu^{-2}}\cdot
\left\vert\frac{y^{2} +py +q}{y-w}\right\vert,\quad
\forall\, x\in U_{\delta}(x_{i,k})\cap 
\sigma_{2\omega T_{\delta}(x_{i,j})}^{-1}\bigl(U_{\delta}(x_{i,k+1})\bigr),
\end{equation}
where variable $y$ and constants $\mu, a, p, q, w$ are
\begin{align}\label{def_y_mu_a_p_q_w}
\nonumber
&y = r_{i,k}^{\max}\bigl(x - x_{i,k}^{(*)}\bigr),\quad
\mu = \left(C_{2} g\right)^{2T_{\delta}(x_{i,k+1})-1}
\left(1 + O\left(g^{-8\bigl(T_{\delta}(x_{i,k})-T_{\delta}(x_{i,k+1})\bigr)}\right)\right),\\
&a = -\frac{r_{i,k+1}^{\max}}{r_{i,k}^{\min}},\quad
p = r_{i,k}^{\min} \Delta_{i,k},\quad
q = \frac{r_{i,k}^{\min}}{r_{i,k+1}^{\max}}g^{-2}\mu^{-2},\quad
w = \frac{r_{i,k+1}^{\max}}{r_{i,k}^{\min}\mu^{2} - r_{i,k+1}^{\max}}.
\end{align}
Here we used the second condition $(\ref{cond_L7})$ to obtain $a<0$.

Finally, we notice that under condition $p^{2} < q$ the minimum value of 
the right hand side in inequality $(\ref{est_Phi_ik_2})$ is $\vert a\vert^{1/2}\mu^{-1}\bigl(1+O(\mu^{-2})\bigr)$. Substituting this value into $(\ref{est_Phi_ik_2})$ and taking into account $(\ref{def_y_mu_a_p_q_w})$ one finishes the proof. $\qed$

It follows from $(\ref{est_m})$ and $(\ref{est_lambda_ik})$ that 
$\forall\, x\in U_{\delta}(x_{i,k})\cap 
\sigma_{2\omega T_{\delta}(x_{i,k})}^{-1}\bigl(U_{\delta}(x_{i,k+1})\bigr)$
\begin{equation*}
\lambda_{i,k}^{(2)}(x) \ge \frac
{\Lambda_{0}\left(C_{1}\delta\right)^{2(T_{\delta}(x_{i,k+1})-1)}}
{\left(C_{2}\right)^{2T_{\delta}(x_{i,k+1})-1}}
g^{2(T_{\delta}(x_{i,k}) - T_{\delta}(x_{i,k+1}))}.
\end{equation*}
If periods $T_{\delta}(x_{i,k})$, $T_{\delta}(x_{i,k+1})$ satisfy
\begin{equation*}\label{est_T_k_k+1}
T_{\delta}(x_{i,k}) \ge 2 T_{\delta}(x_{i,k+1}),
\end{equation*}
then 
\begin{equation*}
\lambda_{i,k}^{(2)}(x)\cdot 
\left\vert \cot\Phi_{i,k}^{(2)}(x)\right\vert \ge C_{\lambda}^{(2)} g,\quad
\forall\, x\in U_{\delta}(x_{i,k})\cap 
\sigma_{2\omega}^{-T_{\delta}(x_{i,k})}\bigl(U_{\delta}(x_{i,k+1})\bigr)
\end{equation*}
with some positive constant $C_{\lambda}^{(2)}$. It means that the product 
$P_{i,k}^{(2)}$ satisfies conditions of Lemma 3. 

{\bf Remark:} It has to be noted that condition $(\ref{est_T_k_k+1})$ is stronger than condition 5 in $(\ref{cond_C2_2})$. The necessity to use such kind of condition is dictated by the fact that Lemma 5 cannot be applied to the product $P_{i,k}^{(2)}$ directly. Indeed, inequality $(\ref{cond_L5_2})$ does not hold for all $x\in U_{\delta}(x_{i,k})\cap 
\sigma_{2\omega}^{-T_{\delta}(x_{i,k})}\bigl(U_{\delta}(x_{i,k+1})\bigr)$. It is satisfied only for those $x$ which are sufficiently close to $x_{i,k}^{(*)}$.

For a fixed $i$ we consider the equivalence class $L_{i}$ of cyclically connected points. 
\begin{defs}
We say that an equivalence class $L_{i}$ obeys a condition $(H_{4})$ if
\begin{equation*}
1.\quad\quad\quad 
b'\left(x_{i,k}^{(0)}\right)\cdot b'\left(x_{i,k+1}^{(0)}\right) < 0,\quad
\forall\, k = 1,\ldots, S_{i}
\end{equation*}
and one of two sets of conditions:  
\begin{eqnarray*}
2.1\quad\quad &a.&\quad\quad 
T_{\delta}(x_{i,2k-1}) \ge 2 T_{\delta}(x_{i,2k}),\quad \forall\, k=1,\ldots, S_{i}/2,\\
&b.&\quad\quad
\vert \Delta_{i,2k-1} \vert \le 
\frac{\left(r_{i,2k-1}^{\max} \cdot r_{i,2k}^{\max}\right)^{-1/2}}
{g(C_{2}g)^{2 T_{\delta}(x_{i,2k})-1}},\quad \forall\, k=1,\ldots, S_{i}/2;\\
\\
2.2\quad\quad &a.&\quad\quad 
T_{\delta}(x_{i,2k}) \ge 2 T_{\delta}(x_{i,2k+1}),\quad \forall\, k=1,\ldots, S_{i}/2,\\
&b.&\quad\quad
\vert \Delta_{i,2k} \vert \le 
\frac{\left(r_{i,2k}^{\max} \cdot r_{i,2k+1}^{\max}\right)^{-1/2}}
{g(C_{2}g)^{2 T_{\delta}(x_{i,2k+1})-1}},\quad \forall\, k=1,\ldots, S_{i}/2;
\end{eqnarray*}
holds.
\end{defs}
{\bf Remarks:} 1. We point out that the first condition of Definition 4 together with assumption $(H'_{3})$ guarantee that $S_{i}$ is even. Thus, conditions 2.1, 2.2 are well-posed.

2. Assume that for a fixed $i\in \{1,\ldots, M\}$ the equivalence class $L_{i}$ obeys condition $(H_{4})$ such that the set of conditions 2.1 of Definition 4 is satisfied. Then one may assert that for all 
$j=1,\ldots, S_{i}/2$ a product of matrices
\begin{equation*}
\vprod\limits_{k=1}^{j} P_{i,2k-1}^{(2)}(x),\quad
x\in \mathcal{M}_{i} = \bigcap\limits_{l=1}^{2j} 
\sigma_{2\omega \tau_{l}}^{-1}
\Bigl(U_{\delta}(x_{i,l})\Bigr),\quad \tau_{l} = \sum\limits_{s=1}^{l-1}T_{\delta}(x_{i,s})
\end{equation*} 
admits estimate $(\ref{coc_est})$. It also has to be noted that set $\mathcal{M}_{i}$ is not empty. Moreover, estimates 2.1b yield
\begin{equation*}
U_{\delta/2}(x_{i,1})\subset \mathcal{M}_{i}\subset U_{\delta}(x_{i,1}).
\end{equation*} 

3. Notice that if $x\notin U_{\delta/2}(\mathcal{C}_{0})$ is such that 
$\sigma_{\omega}(x)\notin U_{\delta/2}(\mathcal{C}_{0})$ then due to $(\ref{cond_b})$ one gets
\begin{equation*}
\lambda(x) \ge \left(\frac{1}{2}C_{1}\delta g\right)^{2},\quad
\left\vert \cot \Phi(x)\right\vert \ge \frac{1}{2}C_{1}\delta g.
\end{equation*}
Hence, if a trajectory of a point $x\in \mathbb{T}^{1}$ under rotation over an angle $-2\pi\omega$ does not fall into $U_{\delta/2}(\mathcal{C}_{0})$ then the corresponding product of matrices $(\ref{coc_B_3})$ satisfies conditions of Lemma 3 provided $g$ to be sufficiently large. 

Summarizing the obtained results, we arrive at a theorem.

\begin{thms}
Let $g_{0}$ be a large positive parameter and $b\in C^{1}(\mathbb{T}^{1},\mathbb{R})$ satisfies condition $(H'_{3})$. Suppose that for $\omega\in (0,1)\setminus\mathbb{Q}$ there exist positive constants $\delta$ and $h_{*}$ such that
\begin{equation*}
(H''_{1,0})\quad\quad
\min\limits_{x\in U_{\delta}\left(\mathcal{C}_{0}^{(0)}\right)}
\left(\left( b(x) - b(x-\omega)\right)^{2} + g_{0}^{2}b^{2}(x)b^{2}(x-\omega)\right) \ge 
h_{*}^{2}.
\end{equation*}
Assume that $g_{0}$ and $\delta$ are such that the estimates $(\ref{cond_b})$ hold true and the map $R_{\delta}$ is correctly defined for all $g\in [g_{0}, +\infty)$. 

Let $M$ denotes the number of equivalence classes of cyclically connected points. If for some $g\in [g_{0}, +\infty)$ and for every $i\in \{1,\ldots, M\}$ an equivalence class $L_{i}$ obeys condition $(H_{4})$, then
a continued fraction 
$\{(-1, gb_{j}), f_{j}\}$ with $b_{j} = b(x-(j-1)\omega)$ converges and 
$f_{*} = \lim\limits_{n\to \infty}f_{n}\in C(\mathbb{T}^{1},\mathbb{\hat R})$.
\end{thms}
{\bf PROOF:} Let $B$ be the linear transformation $(\ref{coc_B_3})$ of a cocycle associated to the continued fraction $\{(-1, gb_{j}), f_{j}\}$. One notes that condition $(H''_{1,0})$ implies condition $(H''_{1})$ for all $g\in [g_{0}, +\infty)$, which guarantees the existence of positive constant $\Lambda_{0}$ such that 
\begin{equation*}
\lambda(x) \ge \Lambda_{0}g,\quad \forall\, x\in\mathbb{T}^{1}.
\end{equation*}

Consider an arbitrary point 
$x\in \mathbb{T}^{1}$. As it was mentioned above, a sequence of matrices associated to a finite trajectory of $x$ before its first enter to 
$U_{\delta}\left(\mathcal{C}_{0}^{(0)}\right)$ satisfies $(\ref{coc_est})$ (due to Lemma 3). As the trajectory enters $U_{\delta}\left(\mathcal{C}_{0}\right)$ it stays in a $\delta$-neighborhood of critical points $x_{i,k}, k=1,\ldots, S_{i}$ for some particular $i\in \{1,\ldots, M\}$. We remark that in this case the trajectory may approach points 
$x_{i,k}^{(*)}, k=1,\ldots, S_{i}$ very close. Nevertheless, if the number of iteration is multiple of $S_{i}$, then the corresponding product of matrices $B$ still admits estimate
$(\ref{coc_est})$ (due to condition $(H_{4})$). On the other hand, if the trajectory stays in 
$\mathbb{T}^{1}\setminus U_{\delta/2}\left(\mathcal{C}_{0}\right)$ it also satisfies
$(\ref{coc_est})$ (due to Lemma 3).

It  has to be emphasized that dynamics under rotation over the angle $-2\pi\omega$ possesses the following property of monotonicity. Denote
\begin{equation*}
\mathcal{C}_{0, i} = \bigcup\limits_{k=1}^{S_{i}}\{x_{i,k}\}.
\end{equation*}
Then, if the trajectory of a point $x$ leaves $U_{\delta/2}\left(\mathcal{C}_{0,i}\right)$, it does not enter this neighborhood before it leaves 
$U_{\delta}\left(\mathcal{C}_{0,i}\right)$ (see e.g.\cite{MelStr}).

Let $K$ be a number of iteration between the trajectory successively enters and exits $U_{\delta}\left(\mathcal{C}_{0,i}\right)$. We represent it as 
\begin{equation*}
K = m_{i}S_{i} + n_{i},\quad m_{i}, n_{i}\in \mathbb{N}.
\end{equation*}
Notice that the last $n_{i}$ iterations the trajectory stays in 
$U_{\delta}\left(\mathcal{C}_{0,i}\right)\setminus 
U_{\delta/2}\left(\mathcal{C}_{0,i}\right)$. This finishes the proof. $\qed$

{\bf Remarks:} 1. We point out that condition $(H''_{1})$ follows from
\begin{equation*}
(H'''_{1})\quad\quad
\min\limits_{x\in U_{\delta}\left(\mathcal{C}_{0}^{(0)}\right)}
\left\vert b(x) - b(x-\omega)\right\vert \ge h_{*}>0,
\end{equation*}
which does not contain the parameter $g$ and can be easily checked.

2. It has to be emphasized that points $x_{i,k}^{(*)}$ for all $i\in [1,\ldots, M], k\in [1,\ldots, S_{i}]$ depend on parameter $g$. This implies the distances 
$dist(x_{i,k}^{(*)}, x_{i,k+1}^{(*)})$ vary with parameter $g$. Taking into account 
Lemma 6, one estimates
\begin{equation*}
\frac{{\rm d}}{{\rm d}g}dist(x_{i,k}, x_{i,k+1}) = O(g^{-3}).
\end{equation*}
On the other hand, condition $(H_{4})$ prescribes relative positions of 
$x_{i,k}^{(*)}, x_{i,k+1}^{(*)}$ with very high precision. In particular, if we denote
\begin{equation*}
\tau_{\max} = \max\limits_{i=1,\ldots, M}\sum\limits_{k=1}^{S_{i}-1}T_{\delta}(x_{i,k}),\quad
r_{\max} = \max\limits_{i,k} r_{i,k}^{\max},
\end{equation*}
where $r_{i,k}^{\max}$ are defined by $(\ref{est_r_ik})$, then
relative positions of $x_{i,k}^{(*)}, x_{i,k+1}^{(*)}$ are defined up to 
\begin{equation*}
\Delta_{\min} = \frac{r_{\max}^{-2}}{g\left(C_{2}g\right)^{2\tau_{\max}/3-1}}.
\end{equation*}
This suggests that whenever condition $(H_{4})$ holds for some $g_{*} \ge g_{0}$ then one may expect that $(H_{4})$ also holds for the nearest values of $g$ such that
\begin{equation*}
\vert g - g_{*}\vert = O\left(g_{*}^{3-2\tau_{\max}/3}\right).
\end{equation*} 

3. The simplest case corresponds to those functions $b\in C^{1}(\mathbb{T}^{1}, \mathbb{R})$, which have only two points $x_{0}^{(0)}, x_{1}^{(0)}$ satisfying 
\begin{equation*}
b(x_{p}^{(0)}) = 0,\quad b'(x_{p}^{(0)})\neq 0.
\end{equation*}
In this case the critical set $\mathcal{C}_{0}$ also consists of two points. Moreover, if the time of secondary collisions less than the time of primary ones, there is only one equivalence class of cyclically connected points ($M=1$). This case has been studied in \cite{Iva23}.

%  citations should be arranged by order of appearance

%  citations should be arranged by order of appearance


\begin {thebibliography}{9}
\bibitem{BuFe94}{V. S. Buslaev, A. A.  Fedotov,  "Monodromization and Harper equation", S\'em. sur les \'Eq. aux D\'er. Part. {\bf XXI}, 23 pp., \'Ecole Polytech., Palaiseau (1994). http://archive.numdam.org/article/SEDP\_1993-1994\_\_\_\_A22\_0.pdf}
\bibitem{BuFe96}{V. S. Buslaev, A. A. Fedotov,  "Bloch solutions for difference equations", St. Petersburg Math. J. {\bf 7} (4), 561 -- 594 (1996).}
\bibitem{Fed}{A. A. Fedotov, "Monodromization method in the theory of almost-periodic equations, St.-Petersburg Math. J. {\bf 25}(2), 303--325 (2014). https://doi.org/10.1090/S1061-0022-2014-01292-7}
\bibitem{ABD}{A. Avila, J. Bochi, D. Damanik, "Opening gaps in the spectrum of strictly ergodic Schr\"dinger operators", J. Eur. Math. Soc. {\bf 14}, 61--106 (2012). https://doi.org/10.4171/JEMS/296}
\bibitem{Russ}{H. R$\ddot{\rm u}$ssmann, "On optimal estimates for the solutions of linear difference equations on the circle", Celestial. Mech. {\bf 14}, 33--37 (1976). https://doi.org/10.1007/BF01247129}
\bibitem{Avi}{A. Avila, Almost reducibility and absolute continuity I,  arXiv: 1006.0704. 2010. 
https://doi.org/10.48550/arXiv.1006.0704}
\bibitem{Per}{O. Perron, \textit{Die Lehre von Den Kettenbr$\ddot{\rm u}$chen} (Verlag von B. G. Teubner, Leipzig 1913).}
\bibitem{Khi}{A. Ya. Khinchin, \textit{Continued fractions} (The University of Chicago Press, Chicago, Ill.-London, 1964).}
\bibitem{JoTh}{W. M. Jones, W. J. Thron, \textit{Continued fractions. Analytic theory and applications} (Addison-Wesley Publ. Co.,  MA 1980).}
\bibitem{Brez}{C. Brezinski, \textit{History of continued fractions and Pad\'e approximants} (Springer-Verlag, Berlin Heidelberg 1991).}
\bibitem{Yoc}{J.-C. Yoccoz, "Some questions and remarks about SL(2,$\mathbb{R}$) cocycles", In \textit{Modern dynamical systems and applications}, pp.\; 447--458, Cambridge Univ. Press, Cambridge. 2004}
\bibitem{AvBo}{A. Avila, J. Bochi, "A uniform dichotomy for generic SL(2,$\mathbb{R}$)-cocycles over a minimal base", Bull. Soc. Math. France {\bf 135}, 407--417 (2007). https://doi.org/10.24033/bsmf.2540}
\bibitem{Oseled}{V. I. Oseledets, "A multiplicative ergodic theorem. Lyapunov characteristic numbers for dynamical systems", Trans. Mosc. Math. Soc. {\bf 19}, 197--231 (1968).}
\bibitem{Her}{M. Herman, "Une methode pour minorer les exposants de Lyapunov et quelques exemples montrant le charactere local d'un theoreme d'Arnold et de Moser sur le tore en dimension 2", Commun. Math. Helv. {\bf 58}, 453--502 (1983). https://doi.org/10.1007/BF02564647}
\bibitem{LSY}{L.-S. Young, "Lyapunov exponents for some quasi-periodic cocycles", Ergod. Th. \& Dynam. Sys. {\bf 17}, 483--504 (1997). https://doi.org/10.1017/S0143385797079170}
\bibitem{Lya}{M. A. Lyalinov, N. Y., Zhu, "A solution procedure for second-order difference  equations and its application to  electromagnetic-wave diffraction  in a wedge-shaped region", Proc. R. Soc. Lond. A {\bf 459}(2040), 3159--3180 (2003). https://doi.org/10.1098/rspa.2003.1165}
\bibitem{Iva21}{A. V. Ivanov, "On singularly perturbed linear cocycles over irrational rotations", Reg. \& Chaotic Dyn. {\bf 26}(3), 205--221 (2021). https://doi.org/10.1134/S1560354721030011}
\bibitem{Iva23}{A. V. Ivanov, "On $SL(2,\mathbb{R})$-cocycles over irrational rotations with secondary collisions", Reg. \& Chaotic Dyn. {\bf 28} (2), 207 -- 226 (2023). https://doi.org/10.1134/S1560354723020053}
\bibitem{MelStr}{W. Melo, S. Strien, \textit{One-Dimensional Dynamics} (A Ser. Modern Surveys in Math., Vol. 25, Springer-Verlag 1993).}
\end{thebibliography}
\end {document}